\documentclass{article}
\usepackage{amssymb}
\usepackage{amsmath}

\newtheorem{theorem}{Theorem}

\newtheorem{remark}[theorem]{Remark}

\newcommand{\slg}{\mbox{\bfseries\slshape g}}
\input{tcilatex}

\begin{document}

\title{ Geometric and Extensor Algebras and the Differential Geometry of
Arbitrary Manifolds. }
\author{{\footnotesize V. V. Fern\'{a}ndez}$^{1},${\footnotesize A. M. Moya}$%
^{2}${\footnotesize and W. A. Rodrigues Jr}.$^{1}${\footnotesize \ } \\
$^{1}\hspace{-0.1cm}${\footnotesize Institute of Mathematics, Statistics and
Scientific Computation}\\
{\footnotesize \ IMECC-UNICAMP CP 6065}\\
{\footnotesize \ 13083-859 Campinas, SP, Brazil }\\
{\footnotesize e-mail: walrod@ime.unicamp.br virginvelfe@accessplus.com.ar}\\
{\footnotesize \ }$^{2}${\footnotesize Department of Mathematics, University
of Antofagasta, Antofagasta, Chile} \\
{\footnotesize e-mail: mmoya@uantof.cl}}
\maketitle

\begin{abstract}
We give in this paper which is the third in a series of four a theory of
covariant derivatives of representatives of multivector and extensor fields
on an arbitrary open set $U\subset M$, based on the geometric and extensor
calculus on an arbitrary smooth manifold $M$. This is done by introducing
the notion of a connection extensor field $\gamma $ defining a parallelism
structure on $U\subset M$, which represents in a well defined way the action
on $U$ of the restriction there of some given connection $\nabla $ defined
on $M$. Also we give a novel and intrinsic presentation (i.e., one that does
not depend on a chosen orthonormal moving frame) of the torsion and
curvature fields of Cartan's theory. Two kinds of Cartan's connection 
\textit{operator} fields are identified, and both appear in the intrinsic
Cartan's structure equations satisfied by the Cartan's torsion and curvature
extensor fields. We introduce moreover a metrical extensor $g$ in $U$
corresponding to the restriction there of given metric tensor 
\slg
defined on $M$ and also introduce the concept a \textit{geometric structure }%
$(U,\gamma ,g)$ for $U\subset $ $M$ and study metric compatibility of
covariant derivatives induced by the connection extensor $\gamma $. This
permits the presentation of the concept of gauge (deformed) derivatives
which satisfy noticeable properties useful in differential geometry and
geometrical theories of the gravitational field. Several derivatives
operators in metric and geometrical structures, like ordinary and covariant
Hodge coderivatives and some duality identities are exhibit.
\end{abstract}

\tableofcontents

\section{Introduction}

This is the third paper in a series of four, where we continue our
exposition of how to use Clifford and extensor algebras methods in the study
of the differential geometry of an of a $n$-dimensional smooth manifold $M$
of arbitrary topology, supporting a metric field $%
\slg%
$ \ (of signature $(p,q)$, $p+q=n$) and an arbitrary connection $\nabla $.
It has three main sections (2, 3 and 4), besides the introduction,
conclusions and an Appendix. Section 2 is dedicated to the description of
the theory of covariant derivatives of representatives of multivector and
extensor fields in our formalism. In Section 2.1 we recall how choosing,
like in \cite{2} a chart $(U_{o},\phi _{o})$ of a given atlas of $M$ and the
associated canonical vector space $\mathcal{U}_{o}$ we can represent the
effects \ of \ the restriction $\left. \nabla \right\vert _{U}$ on $U$ of
any connection on $M$ defining a \emph{parallelism structure} by a \emph{%
connection extensor field} $\gamma $ on $U\subset U_{o}$.

In Sections 2.3 and 2.4 we introduce the concepts of $a$-directional \textit{%
covariant} derivatives (which represents $\nabla _{a}$ on $U\subset U_{o}$)
of \ the representatives of multivector and extensor fields, respectively,
and prove the main properties satisfied by these objects. In Section 2.5 we
give a thoughtful study of the so-called \emph{symmetric parallelism
structures}, where among others we present and prove a \emph{Bianchi-like
identity} and give an intrinsic Cartan theory of the torsion and curvature
extensor fields. \emph{Cartan's connections of the first and second kind }%
are identified, and both appear on our version of Cartan's structure
equations. We emphasize the novelty of our approach to Cartan theory, namely
that it does not depend on any \textit{chosen} orthonormal moving frame,
hence, the name\textit{\ intrinsic} used above. In the Appendix some
examples are worked in detail in order to show how some of the concepts
developed in Section 2 are related to standard ones which deals with the
same subject.

In Section 3 using previous results (\cite{1,2} and Section 2) we first
introduce (Section 3.1) a metric structure on a smooth manifold through the
concept of a metric extensor field $g$ associated to a given metric tensor $%
\slg%
$. Christoffel operators\footnote{%
These objects generalize the Christofell symbols of the standard formalism
that are defined only for vector fields of a coordinate basis.} and the
associated Levi-Civita connection field for \ $U\subset U_{o}\subset M$ are
introduced and their properties are given. The structure of the Levi-Civita
connection field is shown to consist of two pieces, a $g$-symmetric and a $g$%
-skewsymmetric parts which have deep geometric and algebraic meanings. In
Section 3.2 we study metric compatible covariant derivatives. There, the
relationship between the connection extensor fields on $U$ and covariant
derivatives corresponding to \textit{deformed} (metric compatible)
geometrical structures $(U,\gamma ,g)$ are given and analyzed. The crucial
result of this section is Theorem 2 which relates pair of deformed covariant
derivative operators associated to different deformed metrics, and is thus
an important tool for correct approximations in geometric theories of the
gravitational field where some metric is supposed to be a `small'
deformation of some other.

In Section 4 we continue the development of our theory of multivector and
extensor calculus on smooth manifolds, introducing in Section 4.1 the
concept of ordinary Hodge coderivatives, duality identities, and Hodge
coderivative identities. Then in Section 4.2 we analyzed deeply the concept
of a Levi-Civita connection in our formalism and the remarkable concept of 
\textit{gauge }deformed derivatives. Several important formulas that appear
in the Lagrangian formulation of the theory of multivector and extensor
fields on smooth manifolds are obtained.\footnote{%
In \cite{3} we gave a preliminary presentation of the Lagrangian theory of
multivector and extensor fields on Minkowski spacetime. See also related
developments in \cite{rodoliv2006}.} In Section 4.3 we introduce the concept
of \emph{covariant} Hodge coderivative. We study in details how all these
important concepts are related and how they can be utilized to clarify
several issues in geometrical theories of the gravitational field.

We emphasize, as already done in \cite{1} that the methods introduced in
this paper are completely general and applies to the study of the geometry
of any arbitrary smooth manifold of arbitrary topology. And indeed, suppose
that the canonical vector space $\mathcal{U}_{o}$ associated to the chart $%
(U_{o},\phi _{o})$ of a given atlas is not enough to perform calculations
involving a region $V$ outsided $U_{o}$. In this case, all we need to do is
to choose another chart $(U_{1},\phi _{1})$ of the atlas with coordinates $%
\{x^{\prime \mu }\}$ and such that $V\subset U_{1}$ and construct \ a
geometrical algebra associated to the canonical space $\mathcal{U}_{1}$
determined by $(U_{1},\phi _{1})$. Of course, if the manifold is toroidal
then an unique chart may cover it and life will be simpler, in the sense
that some global questions can also be discussed with only the introduction
of the canonical geometrical algebra of a unique canonical space.

\section{Covariant Derivatives of Multivector and Extensor Fields}

\subsection{Representation of a Parallelism Structure on $U$}

Let $M$ be a $n$-dimensional smooth manifold. Then, as already recalled in 
\cite{1} for any point $o\in M$ there exists a local coordinate system $%
(U_{o},\phi _{o})$ with coordinates $\{x^{\mu }\}$.

As in \cite{1}, let $\mathcal{U}_{o}$ be the \emph{canonical vector space}
for $(U_{o},\phi _{o})$, and $U$ be an open subset of $U_{o}.$ We denote the
ring (with identity) of smooth scalar fields on $U,$ the module of
(representatives) of smooth vector fields on $U$ and the module of
(representatives) of smooth multivector fields on $U$ respectively by $%
\mathcal{S}(U),$ $\mathcal{V}(U)$ and $\mathcal{M}(U).$ The set of
(representatives) of smooth $k$-vector fields on $U$ is denoted by $\mathcal{%
M}^{k}(U).$ The module of smooth $k$-extensor fields on $U$ is denoted by $k$%
-$ext(\mathcal{M}_{1}^{\diamond }(U),\ldots ,\mathcal{M}_{k}^{\diamond }(U);%
\mathcal{M}^{\diamond }(U))$.

Any smooth vector elementary $2$-extensor field on $U$ is said to be a \emph{%
connection field} on $U$. As we will see in what follows it represents the
effect of the restriction on $U$ of a connection $\nabla $ defined in $M$. A
general connection field will be denoted by $\gamma ,$ i.e., $\gamma
:U\rightarrow 2$-$ext^{1}(\mathcal{U}_{o})$. The smoothness of such $\gamma $
means that for all $a,b\in \mathcal{V}(U)$ the vector field defined by $U\ni
p\mapsto \gamma _{(p)}(a(p),b(p))\in \mathcal{U}_{o}$ is itself smooth.

The open set $U$ equipped with such a connection field $\gamma ,$ namely $%
(U,\gamma ),$ will be said to be a \emph{parallelism structure} on $U$
representing \textit{there} the action of a given connection $\nabla $
defined in $M$.

\begin{remark}
Please, take notice that, of course, as defined $\gamma $ cannot be extended
to all $M$. However, this does not reduce in any way its theoretical
importance and as a tool for performing easily very sophisticated
calculations.
\end{remark}

Let us take $a\in \mathcal{U}_{o}.$ A smooth $(1,1)$-extensor field on $U,$
namely $\gamma _{a},$ defined as $U\ni p\mapsto \left. \gamma
_{a}\right\vert _{(p)}\in ext_{1}^{1}(\mathcal{U}_{o})$ such that for all $%
b\in \mathcal{U}_{o}$ 
\begin{equation}
\left. \gamma _{a}\right\vert _{(p)}(b)=\gamma _{(p)}(a,b),  \label{PS.1}
\end{equation}%
will be called an $a$-\emph{directional connection field on }$U$, of course,
associated to $(U,\gamma ).$

We emphasize that the $(1,1)$-extensor character and the smoothness of $%
\gamma_{a}$ are immediate consequences of the vector elementary $2$-extensor
character and the smoothness of $\gamma.$

A smooth $(1,2)$-extensor field on $U,$ namely $\Omega ,$ defined as $U\ni
p\mapsto \Omega _{(p)}\in ext_{1}^{2}(\mathcal{U}_{o})$ such that for all $%
a\in \mathcal{U}_{o}$ 
\begin{equation}
\Omega _{(p)}(a)=\dfrac{1}{2}biv[\left. \gamma _{a}\right\vert _{(p)}],
\label{PS.2}
\end{equation}%
will be called (for reason that will become clear in what follows) a \emph{%
gauge connection field on }$U.$

From the definition of $biv[t]$ (see \cite{2}), taking any pair of
reciprocal frame fields on $U,$ say $(\{e_{\mu}\},\{e^{\mu}\}),$ and using
Eq.(\ref{PS.1}), we can write Eq.(\ref{PS.2}) as 
\begin{equation}
\Omega_{(p)}(a)=\frac{1}{2}\gamma_{(p)}(a,e^{\mu}(p))\wedge e_{\mu}(p)=\frac{%
1}{2}\gamma_{(p)}(a,e_{\mu}(p))\wedge e^{\mu}(p).  \label{PS.2a}
\end{equation}
So, we see that the $(1,2)$-extensor character and the smoothness of $\Omega$
are easily deduced from the vector elementary $2$-extensor character and the
smoothness of $\gamma.$

Let us take any pair of reciprocal frame fields on $U,$ say $(\{e_{\mu
}\},\{e^{\mu }\})$, i.e., $e_{\mu }\cdot e^{\nu }=\delta _{\mu }^{\nu }.$
Let $\Gamma _{a}$ be the \emph{generalized }(\emph{extensor field}) of $%
\gamma _{a}$ (see \cite{2}), i.e., $\Gamma _{a}$ defined as $U\ni p\mapsto
\left. \Gamma _{a}\right\vert _{(p)}\in ext(\mathcal{U}_{o}),$ is a smooth
extensor field on $U$ such that for all $X\in \bigwedge \mathcal{U}_{o}$ 
\begin{equation}
\left. \Gamma _{a}\right\vert _{(p)}(X)=\left. \gamma _{a}\right\vert
_{(p)}(e^{\mu }(p))\wedge (e_{\mu }(p)\lrcorner X)=\left. \gamma
_{a}\right\vert _{(p)}(e_{\mu }(p))\wedge (e^{\mu }(p)\lrcorner X).
\label{PS.3}
\end{equation}

It is easily seen that (for a given $X$) the multivector appearing on the
right side of Eq.(\ref{PS.3}) does not depend on the choice of the
reciprocal frame fields. The extensor character and the smoothness of $%
\Gamma _{a}$ follows from the $(1,1)$-extensor character and the smoothness
of $\gamma _{a}.$

We will usually \textit{omit} the letter $p$ in writing the definitions
given by Eq.(\ref{PS.1}), Eq.(\ref{PS.2}) and Eq.(\ref{PS.3}), and other
equations using extensor fields. No confusion should arise with this
standard practice. Eq.(\ref{PS.3}) might be also written in the succinct
form $\Gamma _{a}(X)=\gamma _{a}(\partial _{b})\wedge (b\lrcorner X).$

We end this section by presenting some of the basic properties satisfied by $%
\Gamma_{a}$.\vspace{0.1in}

\textbf{i.} $\Gamma_{a}$ is grade-preserving, i.e., 
\begin{equation}
\text{if }X\in\mathcal{M}^{k}(U),\text{ then }\Gamma_{a}(X)\in\mathcal{M}%
^{k}(U).  \label{PS.4}
\end{equation}

\textbf{ii}. For any $X\in\mathcal{M}(U)$, 
\begin{align}
\Gamma_{a}(\widehat{X}) & =\widehat{\Gamma_{a}(X)},  \label{PS.5a} \\
\Gamma_{a}(\widetilde{X}) & =\widetilde{\Gamma_{a}(X)},  \label{PS.5b} \\
\Gamma_{a}(\overline{X}) & =\overline{\Gamma_{a}(X)}.  \label{PS.5c}
\end{align}

\textbf{iii.} For any $f\in\mathcal{S}(U),$ $b\in\mathcal{V}(U)$ and $X,Y\in%
\mathcal{M}(U),$ 
\begin{align}
\Gamma_{a}(f) & =0,  \label{PS.6a} \\
\Gamma_{a}(b) & =\gamma_{a}(b),  \label{PS.6b} \\
\Gamma_{a}(X\wedge Y) & =\Gamma_{a}(X)\wedge Y+X\wedge\Gamma_{a}(Y).
\label{PS.6c}
\end{align}

\textbf{iv.} The adjoint of $\Gamma_{a},$ namely $\Gamma_{a}^{\dagger},$ is
the generalized of the adjoint of $\gamma_{a}$, namely $\gamma_{a}^{%
\dagger}, $ i.e., 
\begin{equation}
\Gamma_{a}^{\dagger}(X)=\gamma_{a}^{\dagger}(\partial_{b})\wedge(b\lrcorner
X).  \label{PS.7}
\end{equation}

\textbf{v}. The symmetric (skew-symmetric) part of $\Gamma_{a},$ namely $%
\Gamma_{a\pm}=\dfrac{1}{2}(\Gamma_{a}\pm\Gamma_{a}^{\dagger}),$ is the
generalized of the symmetric (skew-symmetric) part of $\gamma_{a},$ namely $%
\gamma_{a\pm}=\dfrac{1}{2}(\gamma_{a}\pm\gamma_{a}^{\dagger}),$ i.e., 
\begin{equation}
\Gamma_{a\pm}(X)=\gamma_{a\pm}(\partial_{b})\wedge(b\lrcorner X).
\label{PS.8}
\end{equation}

\textbf{vi}. $\Gamma_{a-}$ can be factorized by a remarkable formula which
only involves $\Omega.$ It is 
\begin{equation}
\Gamma_{a-}(X)=\Omega(a)\times X.  \label{PS.9}
\end{equation}

\textbf{vii}. For any $X,Y\in\mathcal{M}(U)$ it holds 
\begin{equation}
\Gamma_{a-}(X\ast Y)=\Gamma_{a-}(X)\ast Y+X\ast\Gamma_{a-}(Y),  \label{PS.10}
\end{equation}
where $\ast$ means any suitable product of smooth multivector fields, either 
$(\wedge),$ $(\cdot),$ $(\lrcorner,\llcorner)$ or $(b$-\emph{Clifford product%
}$),$ see \cite{2}.

\subsection{$a$-Directional Covariant Derivatives of Multivector Fields}

Given a parallelism structure $(U,\gamma),$ let us take $a\in\mathcal{U}_{o}$%
. Then, associated to $(U,\gamma)$ we can introduce two $a$-\emph{%
directional covariant derivative operators} ($a$-\emph{DCDO's}), namely $%
\nabla_{a}^{+}$ and $\nabla_{a}^{-},$ which act on the module of smooth
multivector fields on $U.$

They are defined by $\nabla _{a}^{\pm }:\mathcal{M}(U)\rightarrow \mathcal{M}%
(U)$ such that 
\begin{align}
\nabla _{a}^{+}X(p)& =a\cdot \partial _{o}X(p)+\left. \Gamma _{a}\right\vert
_{(p)}(X(p))  \label{CDM.1a} \\
\nabla _{a}^{-}X(p)& =a\cdot \partial _{o}X(p)-\left. \Gamma _{a}^{\dagger
}\right\vert _{(p)}(X(p))\text{ for each }p\in U,  \label{CDM.1b}
\end{align}%
where $a\cdot \partial _{o}$ is the canonical $a$-\emph{DODO} as defined in 
\cite{1}.

We emphasize that each of $\nabla_{a}^{+}$ and $\nabla_{a}^{-}$ satisfies 
\emph{indeed }the fundamental properties which a well-defined covariant
derivative is expected to have. This is trivial to verify whenever we take
into account the well-known properties of $a\cdot\partial_{o}$ (see \cite{1}%
), and the properties of $\Gamma_{a}$ given by Eq.(\ref{PS.4}), Eqs.(\ref%
{PS.6a}), (\ref{PS.6b}) and (\ref{PS.6c}), and Eq.(\ref{PS.7}).

As usual we will write Eq.(\ref{CDM.1a}) and Eq.(\ref{CDM.1b}) by omitting $%
p $ when no confusion arises.

The smooth multivector fields on $U,$ namely $\nabla_{a}^{+}X$ and $\nabla
_{a}^{-}X,$ will be respectively called the \emph{plus }and the \emph{minus} 
$a$-\emph{directional covariant derivatives} of $X.$

We summarize some of the most important properties for the pair of $a$-\emph{%
DCDO's} $\nabla_{a}^{+}$ and $\nabla_{a}^{-}.\vspace{0.1in}$

\textbf{i.} $\nabla_{a}^{+}$ and $\nabla_{a}^{-}$ are grade-preserving
operators on $\mathcal{M}(U),$ i.e., 
\begin{equation}
\text{if }X\in\mathcal{M}^{k}(U),\text{ then }\nabla_{a}^{\pm}X\in \mathcal{M%
}^{k}(U).  \label{CDM.2}
\end{equation}

\textbf{ii. }For all $X\in\mathcal{M}(U),$ and for any $\alpha,\alpha^{%
\prime }\in\mathbb{R}$ and $a,a^{\prime}\in\mathcal{U}_{o}$ we have 
\begin{equation}
\nabla_{\alpha a+\alpha^{\prime}a^{\prime}}^{\pm}X=\alpha\nabla_{a}^{\pm
}X+\alpha^{\prime}\nabla_{a^{\prime}}^{\pm}X.  \label{CDM.3}
\end{equation}

\textbf{iii. }For all $f\in\mathcal{S}(U)$ and $X,Y\in\mathcal{M}(U)$ we
have 
\begin{align}
\nabla_{a}^{\pm}f & =a\cdot\partial_{o}f,  \label{CDM.4a} \\
\nabla_{a}^{\pm}(X+Y) & =\nabla_{a}^{\pm}X+\nabla_{a}^{\pm}Y,  \label{CDM.4b}
\\
\nabla_{a}^{\pm}(fX) & =(a\cdot\partial_{o}f)X+f(\nabla_{a}^{\pm}X).
\label{CDM.4c}
\end{align}

\textbf{iv. }For all $X,Y\in\mathcal{M}(U)$ we have 
\begin{equation}
\nabla_{a}^{\pm}(X\wedge Y)=(\nabla_{a}^{\pm}X)\wedge Y+X\wedge(\nabla
_{a}^{\pm}Y).  \label{CDM.5}
\end{equation}

\textbf{v. }For all $X,Y\in\mathcal{M}(U)$ we have 
\begin{equation}
(\nabla_{a}^{+}X)\cdot Y+X\cdot(\nabla_{a}^{-}Y)=a\cdot\partial_{o}(X\cdot
Y).  \label{CDM.6}
\end{equation}

It should be noticed that $(\nabla_{a}^{+},\nabla_{a}^{-})$ as defined by
Eq.(\ref{CDM.1a}) and Eq.(\ref{CDM.1b}) is the unique pair of $a$-\emph{%
DCDO's} associated to $(U,\gamma)$ which satisfies the remarkable property
given by Eq.(\ref{CDM.6}).

We emphasize here that the $a$-\emph{DODO }$a\cdot\partial_{o}$ acting on $%
\mathcal{M}(U),$ is also a well-defined $a$-\emph{DCDO.} In this particular
case, the connection field $\gamma$ is \emph{identically zero} and the plus
and minus $a$-\emph{DCDO's} are equal to each other, and both of them
coincide with $a\cdot\partial_{o}.$

We introduce yet another well-defined $a$-\emph{DCDO} which acts also on the
module of smooth multivector fields on $U,$ namely $\nabla_{a}^{0}.$

It is defined by 
\begin{equation}
\nabla_{a}^{0}X=\frac{1}{2}(\nabla_{a}^{+}X+\nabla_{a}^{-}X).  \label{CDM.7}
\end{equation}
But, by using Eqs.(\ref{CDM.1a}) and (\ref{CDM.1b}), and Eq.(\ref{PS.9}), we
might write else 
\begin{equation}
\nabla_{a}^{0}X=a\cdot\partial_{o}X+\Omega(a)\times X.  \label{CDM.8}
\end{equation}

The $a$-\emph{DCDO} $\nabla_{a}^{0}$ satisfies the same properties which
hold for each one of the $a$-\emph{DCDO's} $\nabla_{a}^{+}$ and $%
\nabla_{a}^{-}.$ But, it has also an additional remarkable property 
\begin{equation}
(\nabla_{a}^{0}X)\cdot Y+X\cdot(\nabla_{a}^{0}Y)=a\cdot\partial_{o}(X\cdot
Y).  \label{CDM.9}
\end{equation}

Moreover, it satisfies a Leibnitz-like rule for any suitable product of
smooth multivector fields, i.e., 
\begin{equation}
\nabla_{a}^{0}(X*Y)=(\nabla_{a}^{0}X)*Y+X*(\nabla_{a}^{0}Y).  \label{CDM.10}
\end{equation}

\subsubsection{Connection Operators}

Associated to any parallelism structure $(U,\gamma)$ we can introduce two
remarkable operators which map $2$-uples of smooth vector fields to smooth
vector fields.

They are defined by $\Gamma ^{\pm }:\mathcal{V}(U)\times \mathcal{V}%
(U)\rightarrow \mathcal{V}(U)$ such that 
\begin{equation}
\Gamma ^{\pm }(a,b)=\nabla _{a}^{\pm }b,  \label{CO.1}
\end{equation}%
and will be called the \emph{connection operators} of $(U,\gamma )$.

We summarize the basic properties of them.\vspace{0.1in}

\textbf{i.} For all $f\in \mathcal{S}(U),$ and $a,a^{\prime },b,b^{\prime
}\in \mathcal{V}(U),$ we have 
\begin{align}
\Gamma ^{\pm }(a+a^{\prime },b)& =\Gamma ^{\pm }(a,b)+\Gamma ^{\pm
}(a^{\prime },b),  \label{CO.2a} \\
\Gamma ^{\pm }(a,b+b^{\prime })& =\Gamma ^{\pm }(a,b)+\Gamma ^{\pm
}(a,b^{\prime }),  \label{CO.2b} \\
\Gamma ^{\pm }(fa,b)& =f\Gamma ^{\pm }(a,b),  \label{CO.2c} \\
\Gamma ^{\pm }(a,fb)& =(a\cdot \partial _{o}f)b+f\Gamma ^{\pm }(a,b).
\label{CO.2d}
\end{align}%
As we can observe both connection operators satisfy the linearity property
only with respect to the first smooth vector field variable. Thus,
connection operators are \textit{not} extensor fields, because no linear in
the second argument.

\textbf{ii.} For all $a,b,c\in \mathcal{V}(U),$ we have 
\begin{equation}
\Gamma ^{+}(a,b)\cdot c+b\cdot \Gamma ^{-}(a,c)=a\cdot \partial _{o}(b\cdot
c),  \label{CO.3}
\end{equation}%
which is an immediate consequence of Eq.(\ref{CDM.6}).

\subsubsection{Deformation of Covariant Derivatives}

Let $(\nabla _{a}^{+},\nabla _{a}^{-})$ be any pair of $a$-\emph{DCDO's} \
and $\lambda $ a non-singular smooth $(1,1)$-extensor field on $U$ $.$ We
define the deformation of these covariant derivatives as the pair $%
(_{\lambda }\nabla _{a}^{+},_{\lambda }\nabla _{a}^{-})$ by 
\begin{align}
_{\lambda }\nabla _{a}^{+}X& =\underline{\lambda }(\nabla _{a}^{+}\underline{%
\lambda }^{-1}(X)),  \label{CDM.11a} \\
_{\lambda }\nabla _{a}^{-}X& =\underline{\lambda }^{\ast }(\nabla _{a}^{-}%
\underline{\lambda }^{\dagger }(X)),  \label{CDM.11b}
\end{align}%
where $\underline{\lambda }$ is the \emph{extended}\footnote{%
Recall that $\lambda ^{*}=(\lambda ^{-1})^{\dagger }=(\lambda ^{\dagger
})^{-1},$ and $\underline{\lambda }^{-1}=(\underline{\lambda })^{-1}=%
\underline{(\lambda ^{-1})}$ and $\underline{\lambda }^{\dagger }=(%
\underline{\lambda })^{\dagger }=\underline{(\lambda ^{\dagger })},$ see 
\cite{2}$.$} of $\lambda $, is a well-defined pair of $a$-\emph{DCDO's,}
since it satisfies as it is trivial to verify the fundamental properties
given by Eqs.(\ref{CDM.4a}), (\ref{CDM.4b}) and (\ref{CDM.4c}), Eq.(\ref%
{CDM.5}) and Eq.(\ref{CDM.6}). For instance,%
\begin{equation}
_{\lambda }\nabla _{a}^{+}f=\underline{\lambda }(\nabla _{a}^{+}\underline{%
\lambda }^{-1}(f))=a\cdot \partial _{0}f.
\end{equation}%
\ \ 

We verify now that the definitions given by Eq.(\ref{CDM.11a}) and Eq.(\ref%
{CDM.11b}) also satisfy a property analogous to the one given by Eq.(\ref%
{CDM.6}). Indeed, we have 
\begin{align}
(_{\lambda }\nabla _{a}^{+}X)\cdot Y+X\cdot (_{\lambda }\nabla _{a}^{-}Y)&
=(\nabla _{a}^{+}\underline{\lambda }^{-1}(X))\cdot \underline{\lambda }%
^{\dagger }(Y)+\underline{\lambda }^{-1}(X)\cdot (\nabla _{a}^{-}\underline{%
\lambda }^{\dagger }(Y))  \notag \\
& =a\cdot \partial _{o}(\underline{\lambda }^{-1}(X)\cdot \underline{\lambda 
}^{\dagger }(Y)),  \notag \\
& =a\cdot \partial _{o}(X\cdot Y).
\end{align}

\subsection{$a$-Directional Covariant Derivatives of Extensor Fields}

The three $a$-\emph{DCDO's} $\nabla _{a}^{+},$ $\nabla _{a}^{-}$ and $\nabla
_{a}^{0}$ which act on $\mathcal{M}(U)$ can be extended in order to act on
the module of smooth $k$-extensor fields on $U$. For any $t\in k$-$ext(%
\mathcal{M}_{1}^{\diamond }(U),\ldots ,\mathcal{M}_{k}^{\diamond }(U);%
\mathcal{M}^{\diamond }(U))$, we can define exactly $3^{k+1}$ covariant
derivatives, namely $\nabla _{a}^{\sigma _{1}\ldots \sigma _{k}\sigma }t\in
k $-$ext(\mathcal{M}_{1}^{\diamond }(U),\ldots ,\mathcal{M}_{k}^{\diamond
}(U);\mathcal{M}^{\diamond }(U)),$ where each of $\sigma _{1},\ldots ,\sigma
_{k},\sigma $ is being used to denote either $(+),$ $(-)$ or $(0).$ They are
given by the following definition.

For all $X_{1}\in\mathcal{M}_{1}^{\diamond}(U),\ldots,X_{k}\in\mathcal{M}%
_{k}^{\diamond}(U),$ $X\in\mathcal{M}^{\diamond}(U)$ 
\begin{align}
& (\nabla_{a}^{\sigma_{1}\ldots\sigma_{k}\sigma}t)_{(p)}(X_{1}(p),\ldots
,X_{k}(p))\cdot X(p)  \notag \\
& =a\cdot\partial_{o}(t_{(p)}(\ldots)\cdot
X(p))-t_{(p)}(\nabla_{a}^{\sigma_{1}}X_{1}(p),\ldots)\cdot X(p)  \notag \\
& -\cdots-t_{(p)}(\ldots,\nabla_{a}^{\sigma_{k}}X_{k}(p))\cdot
X(p)-t_{(p)}(\ldots)\cdot\nabla_{a}^{\sigma}X(p),  \label{CDE.1}
\end{align}
for each $p\in U$.

As usual when no confusion arises we will write Eq.(\ref{CDE.1}) by omitting 
$p.$

We call the reader's attention that each one of the $\nabla _{a}^{\sigma
_{1}\ldots \sigma _{k}\sigma }t$ defined by Eq.(\ref{CDE.1}) is in fact a
smooth $k$-extensor field. Its $k$-extensor character and smoothness can be
easily deduced from the respective properties of $t.$ We note also that in
the first term on the right side of Eq.(\ref{CDE.1}), $a\cdot \partial _{o}$
refers to the canonical $a$-\emph{DODO} as was defined in \cite{1}.

We notice that any smooth $(1,1)$-extensor field on $U,$ say $t,$ has just $%
3^{1+1}=9$ covariant derivatives. For instance, four important covariant
derivatives of such $t$ are given by 
\begin{align}
(\nabla_{a}^{++}t)(X_{1})\cdot X & =a\cdot\partial_{o}(t(X_{1})\cdot X) 
\notag \\
& -t(\nabla_{a}^{+}X_{1})\cdot X-t(X_{1})\cdot\nabla_{a}^{+}X,
\label{CDE.1a} \\
(\nabla_{a}^{+-}t)(X_{1})\cdot X & =a\cdot\partial_{o}(t(X_{1})\cdot X) 
\notag \\
& -t(\nabla_{a}^{+}X_{1})\cdot X-t(X_{1})\cdot\nabla_{a}^{-}X,
\label{CDE.1b} \\
(\nabla_{a}^{--}t)(X_{1})\cdot X & =a\cdot\partial_{o}(t(X_{1})\cdot X) 
\notag \\
& -t(\nabla_{a}^{-}X_{1})\cdot X-t(X_{1})\cdot\nabla_{a}^{-}X,
\label{CDE.1c} \\
(\nabla_{a}^{-+}t)(X_{1})\cdot X & =a\cdot\partial_{o}(t(X_{1})\cdot X) 
\notag \\
& -t(\nabla_{a}^{-}X_{1})\cdot X-t(X_{1})\cdot\nabla_{a}^{+}X,
\label{CDE.1d}
\end{align}
where $X_{1}\in\mathcal{M}_{1}^{\diamond}(U)$ and $X\in\mathcal{M}^{\diamond
}(U),$

We present now some of the basic properties satisfied by these $a$%
-directional covariant derivatives of smooth $k$-extensor fields.\vspace{%
0.1in}

\textbf{i. }For all $f\in\emph{S}(U),$ and $t,u\in k$-$ext(\mathcal{M}%
_{1}^{\diamond}(U),\ldots,\mathcal{M}_{k}^{\diamond}(U);\mathcal{M}%
^{\diamond }(U)),$ it holds 
\begin{align}
\nabla_{a}^{\sigma_{1}\ldots\sigma_{k}\sigma}(t+u) & =\nabla_{a}^{\sigma
_{1}\ldots\sigma_{k}\sigma}t+\nabla_{a}^{\sigma_{1}\ldots\sigma_{k}\sigma }u,
\label{CDE.2a} \\
\nabla_{a}^{\sigma_{1}\ldots\sigma_{k}\sigma}(ft) & =(a\cdot\partial
_{o}f)t+f(\nabla_{a}^{\sigma_{1}\ldots\sigma_{k}\sigma}t).  \label{CDE.2b}
\end{align}

\textbf{ii.} For all $t\in1$-$ext(\mathcal{M}_{1}^{\diamond}(U);\mathcal{M}%
^{\diamond}(U)),$ it holds a noticeable property 
\begin{equation}
(\nabla_{a}^{\sigma_{1}\sigma}t)^{\dagger}=\nabla_{a}^{\sigma\sigma_{1}}t^{%
\dagger}.  \label{CDE.3}
\end{equation}
Note the inversion between $\sigma_{1}$ and $\sigma$ into the $a$-\emph{%
DCDO's }above. As we can see, the three $a$-\emph{DCDO's }$\nabla_{a}^{++},$ 
$\nabla_{a}^{--}$ and $\nabla_{a}^{00}$ commute indeed with the adjoint
operator $\left. {}\right. ^{\dagger}.$

The proof of the above result is as follows. Let us take $X_{1}\in \mathcal{M%
}_{1}^{\diamond}(U)$ and $X\in\mathcal{M}^{\diamond}(U).$ By recalling the
fundamental property of the adjoint operator \cite{2}, and in accordance
with Eq.(\ref{CDE.1}), we can write 
\begin{align*}
(\nabla_{a}^{\sigma_{1}\sigma}t)^{\dagger}(X)\cdot X_{1} & =(\nabla
_{a}^{\sigma_{1}\sigma}t)(X_{1})\cdot X \\
& =a\cdot\partial_{o}(t(X_{1})\cdot X)-t(\nabla_{a}^{\sigma_{1}}X_{1})\cdot
X-t(X_{1})\cdot\nabla_{a}^{\sigma}X \\
& =a\cdot\partial_{o}(t^{\dagger}(X)\cdot X_{1})-t^{\dagger}(X)\cdot
\nabla_{a}^{\sigma_{1}}X_{1}-t^{\dagger}(\nabla_{a}^{\sigma}X)\cdot X_{1}, \\
& =(\nabla_{a}^{\sigma\sigma_{1}}t^{\dagger})(X)\cdot X_{1}.
\end{align*}
Hence, by non-degeneracy of scalar product, the expected result immediately
follows.

\textbf{iii.} For all $t\in1$-$ext(\mathcal{M}_{1}^{\diamond}(U),\mathcal{M}%
^{\diamond}(U)),$ it holds 
\begin{align}
(\nabla_{a}^{++}t)(X_{1}) & =\nabla_{a}^{-}t(X_{1})-t(\nabla_{a}^{+}X_{1}),
\label{CDE.4a} \\
(\nabla_{a}^{+-}t)(X_{1}) & =\nabla_{a}^{+}t(X_{1})-t(\nabla_{a}^{+}X_{1}),
\label{CDE.4b} \\
(\nabla_{a}^{--}t)(X_{1}) & =\nabla_{a}^{+}t(X_{1})-t(\nabla_{a}^{-}X_{1}),
\label{CDE.4c} \\
\nabla_{a}^{-+}t)(X_{1}) & =\nabla_{a}^{-}t(X_{1})-t(\nabla_{a}^{-}X_{1}).
\label{CDE.4d}
\end{align}

We prove here only Eq.(\ref{CDE.4a}). Let us take $X_{1}\in\mathcal{M}%
_{1}^{\diamond}(U)$ and $X\in\mathcal{M}^{\diamond}(U).$ In accordance with
Eq.(\ref{CDE.1}) and by recalling Eq.(\ref{CDM.6}), we have 
\begin{align*}
(\nabla_{a}^{++}t)(X_{1})\cdot X & =a\cdot\partial_{o}(t(X_{1})\cdot
X)-t(\nabla_{a}^{+}X_{1})\cdot X-t(X_{1})\cdot\nabla_{a}^{+}X \\
& =\nabla_{a}^{-}t(X_{1})\cdot X+t(X_{1})\cdot\nabla_{a}^{+}X \\
& -t(\nabla_{a}^{+}X_{1})\cdot X-t(X_{1})\cdot\nabla_{a}^{+}X \\
& =(\nabla_{a}^{-}t(X_{1})-t(\nabla_{a}^{+}X_{1}))\cdot X.
\end{align*}
Hence, by non-degeneracy of scalar product, it follows what was to be proved.

\subsection{Torsion and Curvature Fields}

Let $(U,\gamma )$ be a parallelism structure on $U$, representing there some
connection $\nabla $ defined in $M$. The \emph{smooth vector elementary} $2$-%
\emph{exform field} on $U,$ namely, $\tau $ such that for all $a,b\in 
\mathcal{V}(U)$ 
\begin{equation}
\tau (a,b)=\nabla _{a}^{+}b-\nabla _{b}^{+}a-[a,b],  \label{TCF.1a}
\end{equation}%
i.e., 
\begin{equation}
\tau (a,b)=\gamma _{a}(b)-\gamma _{b}(a),  \label{TCF.1b}
\end{equation}%
is called the \emph{torsion field} of $(U,\gamma ).$

The \emph{smooth vector elementary} $3$-\emph{extensor field }on $U,$ namely 
$\rho,$ such that for all $a,b,c\in\mathcal{V}(U)$ 
\begin{equation}
\rho(a,b,c)=[\nabla_{a}^{+},\nabla_{b}^{+}]c-\nabla_{[a,b]}^{+}c,
\label{TCF.2a}
\end{equation}
i.e., 
\begin{equation}
\rho(a,b,c)=(a\cdot\partial_{o}\gamma_{b})(c)-(b\cdot\partial_{o}\gamma
_{a})(c)+[\gamma_{a},\gamma_{b}](c)-\gamma_{[a,b]}(c),  \label{TCF.2b}
\end{equation}
will be called the \emph{curvature field} of $(U,\gamma).$

It should be emphasized that the curvature field $\rho $ is \emph{%
skew-symmetric} in the first and second variables, i.e., 
\begin{equation}
\rho (a,b,c)=-\rho (b,a,c).  \label{TCF.3}
\end{equation}

\subsubsection{Symmetric Parallelism Structures}

A parallelism structure $(U,\gamma )$ is said to be \emph{symmetric} if and
only if 
\begin{equation}
\gamma _{a}(b)=\gamma _{b}(a).  \label{SPS.1}
\end{equation}%
As we can easily be verified this condition is completely equivalent to 
\begin{equation}
\nabla _{a}^{+}b-\nabla _{b}^{+}a=[a,b],  \label{SPS.2}
\end{equation}%
for all $a,b\in \mathcal{V}(U)$.

So, taking into account Eq.(\ref{TCF.1a}) and Eq.(\ref{TCF.1b}) we have that
a parallelism structure is symmetric if and only if it is torsionless, i.e., 
\begin{equation}
\tau (a,b)=0.  \label{SPS.3}
\end{equation}

We now present and prove some basic properties of a symmetric parallelism
structure.\vspace{0.1in}

\textbf{i.} The curvature field $\rho$ satisfies the cyclic property 
\begin{equation}
\rho(a,b,c)+\rho(b,c,a)+\rho(c,a,b)=0.  \label{SPS.4}
\end{equation}

The proof is as follows. Recalling Eq.(\ref{TCF.2a}) we can write 
\begin{align}
\rho(a,b,c) &
=\nabla_{a}^{+}\nabla_{b}^{+}c-\nabla_{b}^{+}\nabla_{a}^{+}c-\nabla_{\lbrack
a,b]}^{+}c,  \label{SPS.4a} \\
\rho(b,c,a) &
=\nabla_{b}^{+}\nabla_{c}^{+}a-\nabla_{c}^{+}\nabla_{b}^{+}a-\nabla_{\lbrack
b,c]}^{+}a,  \label{SPS.4b} \\
\rho(c,a,b) &
=\nabla_{c}^{+}\nabla_{a}^{+}b-\nabla_{a}^{+}\nabla_{c}^{+}b-\nabla_{\lbrack
c,a]}^{+}b.  \label{SPS.4c}
\end{align}

By adding Eqs.(\ref{SPS.4a}), (\ref{SPS.4b}) and (\ref{SPS.4c}), wherever by
taking into account Eq.(\ref{SPS.2}), we get 
\begin{align}
& \rho(a,b,c)+\rho(b,c,a)+\rho(c,a,b)  \notag \\
& =\nabla_{a}^{+}(\nabla_{b}^{+}c-\nabla_{c}^{+}b)+\nabla_{b}^{+}(\nabla
_{c}^{+}a-\nabla_{a}^{+}c)+\nabla_{c}^{+}(\nabla_{a}^{+}b-\nabla_{b}^{+}a) 
\notag \\
& -\nabla_{\lbrack a,b]}^{+}c-\nabla_{\lbrack b,c]}^{+}a-\nabla_{\lbrack
c,a]}^{+}b,  \notag \\
& =[a,[b,c]]+[b,[c,a]]+[c,[a,b]].  \label{SPS.4d}
\end{align}
Hence, by recalling the so-called Jacobi's identity for the Lie product of
smooth vector fields \cite{1}, the expected result immediate follows.

\textbf{ii.} The curvature field $\rho$ satisfies the so-called \emph{%
Bianchi's identity}, i.e., 
\begin{equation}
(\nabla_{d}^{+++-}\rho)(a,b,c)+(\nabla_{a}^{+++-}\rho)(b,d,c)+(\nabla
_{b}^{+++-}\rho)(d,a,c)=0.  \label{SPS.5}
\end{equation}

The proof of Eq.(\ref{SPS.5}) is as follows. Let us take $a,b,c,d,w\in 
\mathcal{V}(U)$. Taking into account Eq.(\ref{CDE.1}), and using Eq.(\ref%
{CDM.6}), we have 
\begin{align*}
& (\nabla _{d}^{+++-}\rho )(a,b,c)\cdot w \\
& =d\cdot \partial _{o}(\rho (a,b,c)\cdot w)-\rho (\nabla
_{d}^{+}a,b,c)\cdot w-\rho (a,\nabla _{d}^{+}b,c)\cdot w \\
& -\rho (a,b,\nabla _{d}^{+}c)\cdot w-\rho (a,b,c)\cdot \nabla _{d}^{-}w,
\end{align*}%
i.e., 
\begin{align}
(\nabla _{d}^{+++-}\rho )(a,b,c)& =\nabla _{d}^{+}\rho (a,b,c)-\rho (\nabla
_{d}^{+}a,b,c)  \notag \\
& -\rho (a,\nabla _{d}^{+}b,c)-\rho (a,b,\nabla _{d}^{+}c).  \label{SPS.5a}
\end{align}%
By cycling the letters $a,b,d$ into Eq.(\ref{SPS.5a}), we get 
\begin{align}
(\nabla _{a}^{+++-}\rho )(b,d,c)& =\nabla _{a}^{+}\rho (b,d,c)-\rho (\nabla
_{a}^{+}b,d,c)  \notag \\
& -\rho (b,\nabla _{a}^{+}d,c)-\rho (b,d,\nabla _{a}^{+}c),  \label{SPS.5b}
\\
(\nabla _{b}^{+++-}\rho )(d,a,c)& =\nabla _{b}^{+}\rho (d,a,c)-\rho (\nabla
_{b}^{+}d,a,c)  \notag \\
& -\rho (d,\nabla _{b}^{+}a,c)-\rho (d,a,\nabla _{b}^{+}c).  \label{SPS.5c}
\end{align}

Now, by adding Eqs.(\ref{SPS.5a}), (\ref{SPS.5b}), and (\ref{SPS.5c}),
wherever by using Eq.(\ref{TCF.3}) and Eq.(\ref{SPS.2}), we get 
\begin{align}
& (\nabla_{d}^{+++-}\rho)(a,b,c)+(\nabla_{a}^{+++-}\rho)(b,d,c)+(\nabla
_{b}^{+++-}\rho)(d,a,c)  \notag \\
&
=\nabla_{d}^{+}\rho(a,b,c)+\nabla_{a}^{+}\rho(b,d,c)+\nabla_{b}^{+}%
\rho(d,a,c)  \notag \\
& -\rho([a,b],d,c)-\rho([b,d],a,c)-\rho([d,a],b,c)  \notag \\
& -\rho(a,b,\nabla_{d}^{+}c)-\rho(b,d,\nabla_{a}^{+}c)-\rho(d,a,\nabla
_{b}^{+}c).  \label{SPS.5d}
\end{align}

But, in agreement with Eq.(\ref{TCF.2a}), we can write 
\begin{align}
&
\nabla_{d}^{+}\rho(a,b,c)+\nabla_{a}^{+}\rho(b,d,c)+\nabla_{b}^{+}\rho(d,a,c)
\notag \\
& =[\nabla_{a}^{+},\nabla_{b}^{+}]\nabla_{d}^{+}c+[\nabla_{b}^{+},\nabla
_{d}^{+}]\nabla_{a}^{+}c+[\nabla_{d}^{+},\nabla_{a}^{+}]\nabla_{b}^{+}c 
\notag \\
& -\nabla_{d}^{+}\nabla_{\lbrack a,b]}^{+}c-\nabla_{a}^{+}\nabla_{\lbrack
b,d]}^{+}c-\nabla_{b}^{+}\nabla_{\lbrack d,a]}^{+}c,  \label{SPS.5e}
\end{align}
and 
\begin{align*}
& -\rho([a,b],d,c)-\rho([b,d],a,c)-\rho([d,a],b,c) \\
& =-\nabla_{\lbrack a,b]}^{+}\nabla_{d}^{+}c-\nabla_{\lbrack
b,d]}^{+}\nabla_{a}^{+}c-\nabla_{\lbrack d,a]}^{+}\nabla_{b}^{+}c \\
& +\nabla_{d}^{+}\nabla_{\lbrack a,b]}^{+}c+\nabla_{a}^{+}\nabla_{\lbrack
b,d]}^{+}c+\nabla_{b}^{+}\nabla_{\lbrack d,a]}^{+}c \\
& +\nabla_{\lbrack\lbrack a,b],d]}^{+}c+\nabla_{\lbrack\lbrack
b,d],a]}^{+}c+\nabla_{\lbrack\lbrack d,a],b]}^{+}c.
\end{align*}
i.e., by recalling the Jacobi's identity, 
\begin{align}
& -\rho([a,b],d,c)-\rho([b,d],a,c)-\rho([d,a],b,c)  \notag \\
& =-\nabla_{\lbrack a,b]}^{+}\nabla_{d}^{+}c-\nabla_{\lbrack
b,d]}^{+}\nabla_{a}^{+}c-\nabla_{\lbrack d,a]}^{+}\nabla_{b}^{+}c  \notag \\
& +\nabla_{d}^{+}\nabla_{\lbrack a,b]}^{+}c+\nabla_{a}^{+}\nabla_{\lbrack
b,d]}^{+}c+\nabla_{b}^{+}\nabla_{\lbrack d,a]}^{+}c.  \label{SPS.5f}
\end{align}

Now, by adding Eqs.(\ref{SPS.5e}) and (\ref{SPS.5f}), and using Eq.(\ref%
{TCF.2a}), we get 
\begin{align}
&
\nabla_{d}^{+}\rho(a,b,c)+\nabla_{a}^{+}\rho(b,d,c)+\nabla_{b}^{+}\rho(d,a,c)
\notag \\
& -\rho([a,b],d,c)-\rho([b,d],a,c)-\rho([d,a],b,c)  \notag \\
& =\rho(a,b,\nabla_{d}^{+}c)+\rho(b,d,\nabla_{a}^{+}c)+\rho(d,a,\nabla
_{b}^{+}c).  \label{SPS.5g}
\end{align}

Finally, putting Eq.(\ref{SPS.5g}) into Eq.(\ref{SPS.5d}), the expected
result follows.

\subsubsection{Cartan Fields}

The \emph{smooth} $(1,2)$-\emph{extensor field} on $U,$ namely $\mathbf{%
\Theta }$, defined by 
\begin{equation}
\mathbf{\Theta}(c)=\frac{1}{2}\partial_{a}\wedge\partial_{b}\tau(a,b)\cdot c
\label{CF.1}
\end{equation}
will be called the \emph{Cartan torsion field} of $(U,\gamma).$

We should notice that such $\mathbf{\Theta }$ contains all of the geometric
information which is just contained in $\tau .$ Indeed, Eq.(\ref{CF.1}) can
be inverted in such a way that given any $\mathbf{\Theta }$, there is an
unique $\tau $ that verifies Eq.(\ref{CF.1}). We have, indeed that 
\begin{equation}
\tau (a,b)=\partial _{c}(a\wedge b)\cdot \mathbf{\Theta }(c).  \label{CF.1a}
\end{equation}

The \emph{smooth bivector elementary} $2$-\emph{extensor field }on $U$,
namely $\mathbf{\Omega},$ which is defined by 
\begin{equation}
\mathbf{\Omega}(c,d)=\frac{1}{2}\partial_{a}\wedge\partial_{b}\rho(a,b,c)%
\cdot d  \label{CF.2}
\end{equation}
will be called the \emph{Cartan curvature field }of $(U,\gamma).$

Since Eq.(\ref{CF.2}) can be inverted, by giving $\rho$ in terms of $\mathbf{%
\Omega},$ we see that such $\mathbf{\Omega}$ contains the same geometric
information as $\rho.$ The inversion is realized by 
\begin{equation}
\rho(a,b,c)=\partial_{d}(a\wedge b)\cdot\mathbf{\Omega}(c,d).  \label{CF.2a}
\end{equation}

\subsubsection{Cartan's Structure Equations}

Associated to any parallelism structure $(U,\gamma )$ we can introduce two
noticeable \textit{operators} which map $2$-uples of smooth vector fields to
smooth vector fields. They are:

(a) The mapping $\gamma^{+}:\mathcal{V}(U)\times\mathcal{V}(U)\rightarrow 
\mathcal{V}(U)$ defined by 
\begin{equation}
\gamma^{+}(b,c)=\partial_{a}(\nabla_{a}^{+}b)\cdot c  \label{CSE.1}
\end{equation}
which will be called the \emph{Cartan connection operator of first kind} of $%
(U,\gamma).$

(b) The mapping $\gamma^{-}:\mathcal{V}(U)\times\mathcal{V}(U)\rightarrow 
\mathcal{V}(U)$ defined by 
\begin{equation}
\gamma^{-}(b,c)=\partial_{a}b\cdot(\nabla_{a}^{-}c)  \label{CSE.2}
\end{equation}
which will be called the \emph{Cartan connection operator of second kind} of 
$(U,\gamma).$

We summarize some of the basic properties which are satisfied by the Cartan
operators.\vspace{0.1in}

\textbf{i.} For all $f\in\mathcal{S}(U),$ and $b,b^{\prime},c,c^{\prime}\in%
\mathcal{V}(U)$, we have 
\begin{align}
\gamma^{+}(b+b^{\prime},c) & =\gamma^{+}(b,c)+\gamma^{+}(b^{\prime },c),
\label{CSE.3a} \\
\gamma^{+}(b,c+c^{\prime}) & =\gamma^{+}(b,c)+\gamma^{+}(b,c^{\prime }).
\label{CSE.3b} \\
\gamma^{+}(fb,c) & =(\partial_{o}f)b\cdot c+f\gamma^{+}(b,c),  \label{CSE.3c}
\\
\gamma^{+}(b,fc) & =f\gamma^{+}(b,c).  \label{CSE.3d}
\end{align}

\textbf{ii. }For all $f\in\mathcal{S}(U),$ and $b,b^{\prime},c,c^{\prime}\in%
\mathcal{V}(U)$, we have 
\begin{align}
\gamma^{-}(b+b^{\prime},c) & =\gamma^{-}(b,c)+\gamma^{-}(b^{\prime },c),
\label{CSE.4a} \\
\gamma^{-}(b,c+c^{\prime}) & =\gamma^{-}(b,c)+\gamma^{-}(b,c^{\prime }).
\label{CSE.4b} \\
\gamma^{-}(fb,c) & =f\gamma^{-}(b,c),  \label{CSE.4c} \\
\gamma^{-}(b,fc) & =(\partial_{o}f)b\cdot c+f\gamma^{-}(b,c).  \label{CSE.4d}
\end{align}

We have that the Cartan operator of first kind has the linearity property
with respect to the second variable, and the Cartan operator of second kind
is linear with respect to the first variable.

\textbf{iii. }For any $a,b\in\mathcal{V}(U)$, 
\begin{equation}
\gamma^{+}(b,c)+\gamma^{-}(b,c)=\partial_{o}(b\cdot c).  \label{CSE.5}
\end{equation}

It is an immediate consequence of Eq.(\ref{CDM.6}).\vspace{0.1in}

\textbf{First Cartan's Structure Equation}

For any $c\in \mathcal{V}(U)$ it holds 
\begin{equation}
\mathbf{\Theta }(c)=\partial _{o}\wedge c+\partial _{s}\wedge \gamma
^{-}(s,c),  \label{FCE.1}
\end{equation}%
where $\partial _{o}$ is the Hestenes derivative operator introduced in \cite%
{1}.

We prove Eq.(\ref{FCE.1}) as follows. Using Eq.(\ref{TCF.1a}) we can write 
\begin{align}
\mathbf{\Theta}(c) & =\frac{1}{2}\partial_{a}\wedge\partial_{b}(\nabla
_{a}^{+}b-\nabla_{b}^{+}a-[a,b])\cdot c,  \notag \\
& =\partial_{a}\wedge\partial_{b}(\nabla_{a}^{+}b-a\cdot\partial_{o}b)\cdot
c.  \label{FCE.1a}
\end{align}

A straightforward calculation then yields 
\begin{align}
\partial_{a}\wedge\partial_{b}(\nabla_{a}^{+}b)\cdot c &
=\partial_{a}\wedge\partial_{b}\gamma^{+}(b,c)\cdot a  \notag \\
& =\partial_{b}\wedge\partial_{a}(\gamma^{-}(b,c)\cdot a-a\cdot\partial
_{o}(b\cdot c)),  \notag \\
& =\partial_{b}\wedge\gamma^{-}(b,c)-\partial_{b}\wedge\partial_{o}(b\cdot
c),  \label{FCE.1b}
\end{align}
and 
\begin{align}
-\partial_{a}\wedge\partial_{b}(a\cdot\partial_{o}b)\cdot c & =\partial
_{a}\wedge\partial_{b}(b\cdot(a\cdot\partial_{o}c)-a\cdot\partial_{o}(b\cdot
c))  \notag \\
& =\partial_{a}\wedge(a\cdot\partial_{o}c)+\partial_{b}\wedge\partial
_{a}a\cdot\partial_{o}(b\cdot c),  \notag \\
& =\partial_{o}\wedge c+\partial_{b}\wedge\partial_{o}(b\cdot c).
\label{FCE.1c}
\end{align}

Thus, by putting Eq.(\ref{FCE.1b}) and Eq.(\ref{FCE.1c}) into Eq.(\ref%
{FCE.1a}), we get the expected result.\vspace{0.1in}

\textbf{Second Cartan's Structure Equation}

For any $c,d\in\mathcal{V}(U)$ it holds 
\begin{equation}
\mathbf{\Omega}(c,d)=\partial_{o}\wedge\gamma^{+}(c,d)+\gamma^{+}(c,%
\partial_{s})\wedge\gamma^{-}(s,d).  \label{SCE.1}
\end{equation}

To prove Eq.(\ref{SCE.1})we use Eq.(\ref{TCF.2a}) and write 
\begin{align}
\mathbf{\Omega}(c,d) & =\frac{1}{2}\partial_{a}\wedge\partial_{b}([%
\nabla_{a}^{+},\nabla_{b}^{+}]c-\nabla_{\lbrack a,b]}^{+}c)\cdot d,  \notag
\\
&
=\partial_{a}\wedge\partial_{b}(\nabla_{a}^{+}(\nabla_{b}^{+}c)-\nabla_{a%
\cdot\partial_{o}b}^{+}c)\cdot d.  \label{SCE.1a}
\end{align}

But, by taking a pair of reciprocal frame fields $(\{e_{\sigma}\},\{e^{%
\sigma }\})$ we can write 
\begin{align}
\nabla_{a}^{+}(\nabla_{b}^{+}c)\cdot d & =\nabla_{a}^{+}(\gamma
^{+}(c,e^{\sigma})\cdot be_{\sigma})\cdot d  \notag \\
& =a\cdot\partial_{o}(\gamma^{+}(c,e^{\sigma})\cdot b)e_{\sigma}\cdot
d+\gamma^{+}(c,e^{\sigma})\cdot b\nabla_{a}^{+}e_{\sigma}\cdot d  \notag \\
& =a\cdot\partial_{o}(\gamma^{+}(c,e^{\sigma})\cdot b)e_{\sigma}\cdot
d+\gamma^{+}(c,e^{\sigma})\cdot b\gamma^{+}(e_{\sigma},d)\cdot a,  \notag \\
& =a\cdot\partial_{o}\gamma^{+}(c,e^{\sigma})\cdot be_{\sigma}\cdot
d+\gamma^{+}(c,d)\cdot(a\cdot\partial_{o}b)  \notag \\
& +\gamma^{+}(c,e^{\sigma})\cdot b\gamma^{+}(e_{\sigma},d)\cdot a.
\label{SCE.1b}
\end{align}

Now, the first term into Eq.(\ref{SCE.1a}), by using Eq.(\ref{SCE.1b}) and
the well-known identity $\partial_{o}\wedge(fX)=(\partial_{o}f)\wedge
X+f\partial_{o}\wedge X,$ where $f\in\mathcal{S}(U)$ and $X\in\mathcal{M}%
(U), $ can be written 
\begin{align}
\partial_{a}\wedge\partial_{b}\nabla_{a}^{+}(\nabla_{b}^{+}c)\cdot d &
=\partial_{o}\wedge\gamma^{+}(c,e^{\sigma})(e_{\sigma}\cdot d)+\partial
_{a}\wedge\partial_{b}\gamma^{+}(c,d)\cdot(a\cdot\partial_{o}b)  \notag \\
& -\gamma^{+}(c,e^{\sigma})\wedge\gamma^{+}(e_{\sigma},d),  \notag \\
& =\partial_{o}\wedge\gamma^{+}(c,e^{\sigma})(e_{\sigma}\cdot d)  \notag \\
& -\gamma^{+}(c,e^{\sigma})\wedge(\partial_{o}(e_{\sigma}\cdot d)-\gamma
^{-}(e_{\sigma},d))  \notag \\
& +\partial_{a}\wedge\partial_{b}\gamma^{+}(c,d)\cdot(a\cdot\partial _{o}b),
\notag \\
& =\partial_{o}\wedge\gamma^{+}(c,d)  \notag \\
& +\gamma^{+}(c,e^{\sigma})\wedge\gamma^{-}(e_{\sigma},d)  \notag \\
& +\partial_{a}\wedge\partial_{b}\gamma^{+}(c,d)\cdot(a\cdot\partial_{o}b).
\label{SCE.1c}
\end{align}

It is also 
\begin{equation}
-(\nabla_{a\cdot\partial_{o}b}^{+}c)\cdot d=-\gamma^{+}(c,d)\cdot
(a\cdot\partial_{o}b).  \label{SCE.1d}
\end{equation}

Finally, by putting Eq.(\ref{SCE.1c}) and Eq.(\ref{SCE.1d}) into Eq.(\ref%
{SCE.1a}), we get the expected result.

\subsection{Metric Structure}

Let$\ U$ be an open subset\footnote{%
In this paper we will use the nomenclature and notations just used in \cite%
{1},\cite{2}.} of $U_{o}$. Any symmetric and non-degenerate smooth $(1,1)$%
-extensor field on $U,$ namely $g,$ will be said to be a \emph{metric field}
on $U$. It is quite obvious that it is the extensor representative on $U$ of
a given admissible metric tensor $%
\slg%
$ defined in $M$. This means that $g:U_{o}\rightarrow ext_{1}^{1}(\mathcal{U}%
_{o})$ satisfies $g_{(p)}=g_{(p)}^{\dagger }$ and $\det [g]\neq 0,$ for each 
$p\in U,$ and for all $v\in \mathcal{V}(U)$ the vector field defined by $%
U_{o}\ni p\mapsto g_{(p)}(v(p))$ belongs to $\mathcal{V}(U),$ see \cite{1}.

The open set $U$ equipped with such a metric field $g$, namely $(U,g)$, will
be said to be a \emph{metric structure} on $U$.

The existence of a metric field on $U$ makes possible the introduction of
three kinds of \emph{metric products} of smooth multivector fields on $U$.
These are: (a) the $g$-\emph{scalar product }of $X,Y\in\mathcal{M}(U),$
namely $X\underset{g}{\cdot}Y\in\mathcal{S}(U)$; (b) the \emph{left} and 
\emph{right} $g$-\emph{contracted products }of $X,Y\in\mathcal{M}(U),$
namely $X\underset{g}{\lrcorner}Y\in\mathcal{M}(U)$ and $X\underset{g}{%
\llcorner}Y\in\mathcal{M}(U).$ These products are defined by 
\begin{align}
(X\underset{g}{\cdot}Y)(p) & =\underline{g}_{(p)}(X(p))\cdot Y(p)
\label{MS.1a} \\
(X\underset{g}{\lrcorner}Y)(p) & =\underline{g}_{(p)}(X(p))\lrcorner Y(p)
\label{MS.1b} \\
(X\underset{g}{\llcorner}Y)(p) & =X(p)\llcorner\underline{g}_{(p)}(Y(p)),%
\text{ for each }p\in U.  \label{MS.1c}
\end{align}
Note that in the above formulas $\underline{g}$ is the \emph{\ extended}
(extensor field) of $g$ \cite{2}.

We recall that the $g$-\emph{Clifford product }of $X,Y\in \mathcal{M}(U),$
namely $X\underset{g}{}Y\in \mathcal{M}(U),$ is defined by the following
axioms.

For all $f\in\mathcal{S}(U),$ $b\in\mathcal{V}(U)$ and $X,Y,Z\in \mathcal{M}%
(U)$ 
\begin{align}
f\underset{g}{}X & =X\underset{g}{}f=fX\text{ (scalar multiplication on }%
\mathcal{M}(U)\text{).}  \label{MS.2a} \\
b\underset{g}{}X & =b\underset{g}{\lrcorner}X+b\wedge X,  \label{MS.2b} \\
X\underset{g}{}b & =X\underset{g}{\llcorner}b+X\wedge b.  \label{MS.2c} \\
X\underset{g}{}(Y\underset{g}{}Z) & =(X\underset{g}{}Y)\underset{g}{}Z.
\label{MS.2d}
\end{align}

$\mathcal{M}(U)$ equipped with each one of the products $(\underset{g}{%
\lrcorner})$ or $(\underset{g}{\llcorner})$ is a non-associative algebra
induced by the respective $b$-interior algebra of multivectors. They are
called the $g$-\emph{interior algebras of smooth multivector fields.}

$\mathcal{M}(U)$ equipped with $(\underset{g}{})$ is an associative algebra
fundamentally induced by the $b$-Clifford algebra of multivectors. It is
called the $g$-\emph{Clifford algebra of smooth multivector fields.}

\subsection{Christoffel Operators}

Given a metric structure $(U,g)$ we can introduce the following two
operators which map $3$-uples of smooth vector fields into smooth scalar
fields.

(a) The mapping $\left[ \left. {}\right. ,\left. {}\right. ,\left. {}\right. %
\right] :\mathcal{V}(U)\times \mathcal{V}(U)\times \mathcal{V}(U)\rightarrow 
\mathcal{S}(U)$ defined by 
\begin{align}
\left[ a,b,c\right] & =\frac{1}{2}(a\cdot \partial _{o}(b\underset{g}{\cdot }%
c)+b\cdot \partial _{o}(c\underset{g}{\cdot }a)-c\cdot \partial _{o}(a%
\underset{g}{\cdot }b)  \notag \\
& +c\underset{g}{\cdot }[a,b]+b\underset{g}{\cdot }[c,a]-a\underset{g}{\cdot 
}[b,c])  \label{CHO.1}
\end{align}%
which is called the \emph{Christoffel operator of first kind.} Of course, it
is associated to $(U,g).$

(b) The mapping $\left\{ 
\begin{array}{c}
\\ 
\end{array}%
\right\} :\mathcal{V}(U)\times \mathcal{V}(U)\times \mathcal{V}%
(U)\rightarrow \mathcal{S}(U)$ defined by 
\begin{equation}
\left\{ 
\begin{array}{c}
c \\ 
a,b%
\end{array}%
\right\} =[a,b,g^{-1}(c)]  \label{CHO.2}
\end{equation}%
which is called the\textit{\ }\emph{Christoffel operator of second kind.}

\begin{remark}
Before we proceed we remark that the Christofell operators of first and
second kinds generalize the well known Christofell symbols of the standard
formalism used in the differential geometry that is defined only for the
vector fields of a coordinate basis.
\end{remark}

We summarize now some of the basic properties of the Christoffel operator of
first kind.\vspace{0.1in}

\textbf{i.} For all $f\in \mathcal{S}(U),$ and $a,a^{\prime },b,b^{\prime
},c,c^{\prime }\in \mathcal{V}(U)$, 
\begin{align}
\lbrack a+a^{\prime },b,c]& =[a,b,c]+[a^{\prime },b,c],  \label{CHO.3a} \\
\lbrack fa,b,c]& =f[a,b,c].  \label{CHO.3b} \\
\lbrack a,b+b^{\prime },c]& =[a,b,c]+[a,b^{\prime },c],  \label{CHO.3c} \\
\lbrack a,fb,c]& =f[a,b,c]+(a\cdot \partial _{o}f)b\underset{g}{\cdot }c.
\label{CHO.3d} \\
\lbrack a,b,c+c^{\prime }]& =[a,b,c]+[a,b,c^{\prime }],  \label{CHO.3e} \\
\lbrack a,b,fc]& =f[a,b,c].  \label{CHO.3f}
\end{align}%
These equations say that the Christoffel operator of the first kind has the
linearity property which respect to the first and third smooth vector field
variables.

\textbf{ii.} For all $a,b,c\in\mathcal{V}(U)$, 
\begin{align}
\lbrack a,b,c]+[b,a,c] & =a\cdot\partial_{o}(b\underset{g}{\cdot}%
c)+b\cdot\partial_{o}(c\underset{g}{\cdot}a)-c\cdot\partial_{o}(a\underset{g}%
{\cdot}b)  \notag \\
& +b\underset{g}{\cdot}[c,a]-a\underset{g}{\cdot}[b,c],  \label{CHO.4a} \\
\lbrack a,b,c]-[b,a,c] & =c\underset{g}{\cdot}[a,b].  \label{CHO.4b} \\
\lbrack a,b,c]+[a,c,b] & =a\cdot\partial_{o}(b\underset{g}{\cdot }c),
\label{CHO.4c} \\
\lbrack a,b,c]-[a,c,b] & =b\cdot\partial_{o}(c\underset{g}{\cdot}%
a)-c\cdot\partial_{o}(a\underset{g}{\cdot}b)  \notag \\
& +c\underset{g}{\cdot}[a,b]+b\underset{g}{\cdot}[c,a]-a\underset{g}{\cdot }%
[b,c].  \label{CHO.4d} \\
\lbrack a,b,c]+[c,b,a] & =b\cdot\partial_{o}(c\underset{g}{\cdot }a)+c%
\underset{g}{\cdot}[a,b]-a\underset{g}{\cdot}[b,c],  \label{CHO.4e} \\
\lbrack a,b,c]-[c,b,a] & =a\cdot\partial_{o}(b\underset{g}{\cdot}%
c)-c\cdot\partial_{o}(a\underset{g}{\cdot}b)+b\underset{g}{\cdot}[c,a].
\label{CHO.4f}
\end{align}

\subsection{Levi-Civita Connection Field}

We now introduce a remarkable decomposition of the Christoffel operator of
first kind.

\textbf{Proposition}. There exists a smooth $(1,2)$-extensor field on $U,$
namely $\omega _{0},$ such that the Christoffel operator of first kind can
be written as 
\begin{equation}
\lbrack a,b,c]=\left( a\cdot \partial _{o}b+\frac{1}{2}g^{-1}\circ (a\cdot
\partial _{o}g)(b)+\omega _{0}(a)\underset{g}{\times }b\right) \underset{g}{%
\cdot }c.  \label{LCC.1}
\end{equation}%
Such $\omega _{0}$ is given by 
\begin{equation}
\omega _{0}(a)=-\frac{1}{4}\underline{g}^{-1}(\partial _{b}\wedge \partial
_{c})a\cdot \left( (b\cdot \partial _{o}g)(c)-(c\cdot \partial
_{o}g)(b)\right) .  \label{LCC.2}
\end{equation}

\textbf{Proof}

By using $a\cdot\partial_{o}(X\underset{g}{\cdot}Y)=(a\cdot\partial _{o}X)%
\underset{g}{\cdot}Y+X\underset{g}{\cdot}(a\cdot\partial_{o}Y)+(a\cdot%
\partial_{o}\underline{g})(X)\cdot Y,$ for all $X,Y\in \mathcal{M}(U),$ we
have 
\begin{equation}
a\cdot\partial_{o}(b\underset{g}{\cdot}c)=(a\cdot\partial_{o}b)\underset{g}{%
\cdot}c+b\underset{g}{\cdot}(a\cdot\partial_{o}c)+(a\cdot\partial
_{o}g)(b)\cdot c,  \label{LCC.N1}
\end{equation}
and, by cycling the letters $a,b$ and $c$, we get 
\begin{align}
b\cdot\partial_{o}(c\underset{g}{\cdot}a) & =(b\cdot\partial_{o}c)\underset{g%
}{\cdot}a+c\underset{g}{\cdot}(b\cdot\partial_{o}a)+(b\cdot
\partial_{o}g)(c)\cdot a,  \label{LCC.N2} \\
-c\cdot\partial_{o}(a\underset{g}{\cdot}b) & =-(c\cdot\partial _{o}a)%
\underset{g}{\cdot}b-a\underset{g}{\cdot}(c\cdot\partial_{o}b)-(c\cdot%
\partial_{o}g)(a)\cdot b.  \label{LCC.N3}
\end{align}

A straightforward calculation yields 
\begin{align}
c\underset{g}{\cdot}[a,b] & =c\underset{g}{\cdot}(a\cdot\partial _{o}b)-c%
\underset{g}{\cdot}(b\cdot\partial_{o}a),  \label{LCC.N4} \\
b\underset{g}{\cdot}[c,a] & =b\underset{g}{\cdot}(c\cdot\partial _{o}a)-b%
\underset{g}{\cdot}(a\cdot\partial_{o}c),  \label{LCC.N5} \\
-a\underset{g}{\cdot}[b,c] & =-a\underset{g}{\cdot}(b\cdot\partial _{o}c)+a%
\underset{g}{\cdot}(c\cdot\partial_{o}b).  \label{LCC.N6}
\end{align}

Now, by adding Eqs.(\ref{LCC.N1}), (\ref{LCC.N2}), (\ref{LCC.N3}) and Eqs.(%
\ref{LCC.N4}), (\ref{LCC.N5}), (\ref{LCC.N6}) we get 
\begin{equation*}
2[a,b,c]=2(a\cdot\partial_{o}b)\underset{g}{\cdot}c+(a\cdot\partial
_{o}g)(b)\cdot c+(b\cdot\partial_{o}g)(c)\cdot
a-(c\cdot\partial_{o}g)(a)\cdot b,
\end{equation*}
hence, by taking into account the symmetry property $(n\cdot\partial
_{o}g)^{\dagger}=n\cdot\partial_{o}g,$ it follows 
\begin{equation}
\lbrack a,b,c]=(a\cdot\partial_{o}b)\underset{g}{\cdot}c+\frac{1}{2}%
g^{-1}\circ(a\cdot\partial_{o}g)(b)\underset{g}{\cdot}c+\frac{1}{2}%
a\cdot((b\cdot\partial_{o}g)(c)-(c\cdot\partial_{o}g)(b)).  \label{LCC.N7}
\end{equation}

On another side, a straightforward calculation yields 
\begin{align*}
\omega_{0}(a)\underset{g}{\times}b & =-g(b)\lrcorner\omega_{0}(a) \\
& =\frac{1}{4}g(b)\lrcorner\underline{g}^{-1}(\partial_{p}\wedge\partial
_{q})a\cdot((p\cdot\partial_{o}g)(q)-(q\cdot\partial_{o}g)(p)) \\
& =\frac{1}{4}(b\cdot\partial_{p}g^{-1}(\partial_{q})-b\cdot\partial
_{q}g^{-1}(\partial_{p}))a\cdot((p\cdot\partial_{o}g)(q)-(q\cdot\partial
_{o}g)(p)) \\
& =\frac{1}{2}b\cdot\partial_{p}g^{-1}(\partial_{q})a\cdot((p\cdot
\partial_{o}g)(q)-(q\cdot\partial_{o}g)(p)) \\
& =\frac{1}{2}g^{-1}(\partial_{q})a\cdot((b\cdot\partial_{o}g)(q)-(q\cdot
\partial_{o}g)(b)),
\end{align*}
hence, it follows 
\begin{equation}
(\omega_{0}(a)\underset{g}{\times}b)\underset{g}{\cdot}c=\frac{1}{2}%
a\cdot((b\cdot\partial_{o}g)(c)-(c\cdot\partial_{o}g)(b)).  \label{LCC.N8}
\end{equation}

Finally, putting Eq.(\ref{LCC.N8}) into Eq.(\ref{LCC.N7}), we get the
required result.$\blacksquare$

The \emph{smooth vector elementary} $2$-\emph{extensor field} on $U,$ namely 
$\lambda ,$ defined by 
\begin{equation}
\lambda (a,b)=\frac{1}{2}g^{-1}\circ (a\cdot \partial _{o}g)(b)+\omega
_{0}(a)\underset{g}{\times }b  \label{LCC.3}
\end{equation}%
is a well-defined connection field on $U.$ It will be called the \emph{%
Levi-Civita connection field} on $U$. The open set $U$ endowed with $\lambda 
$, namely $(U,\lambda ),$ will be said to be a \emph{Levi-Civita parallelism
structure} on $U$. It is clear that a \textit{particular } $(U,\lambda )$
which is $g$-antisymmetic may in our formalism be the description in $U$ of
the restriction on $U$ of the usual Levi-Civita connection $D$ of $%
\slg%
$ on $M$.

The $a$-\emph{DCDO's} associated to $(U,\lambda )$, namely $D_{a}^{+}$ and $%
D_{a}^{-}$, are said to be \emph{Levi-Civita }$a$-\emph{DCDO's.} They are
fundamentally defined by $D_{a}^{\pm }:\mathcal{M}(U)\rightarrow \mathcal{M}%
(U)$ such that 
\begin{align}
D_{a}^{+}X& =a\cdot \partial _{o}X+\Lambda _{a}(X),  \label{LCC.3a1} \\
D_{a}^{-}X& =a\cdot \partial _{o}X-\Lambda _{a}^{\dagger }(X),
\label{LCC.3a2}
\end{align}%
where $\Lambda _{a}$ is the \emph{generalized}\footnote{%
We recall from \cite{3} that the generalized of $t\in ext_{1}^{1}(\mathcal{U}%
_{o})$ is $T\in ext(\mathcal{U}_{o})$ given by $T(X)=t(\partial _{n})\wedge
(n\lrcorner X).$} (extensor field) of $\lambda _{a}$.

It should be noted that such a $a$-\emph{DCDO} $D_{a}^{+}$ satisfies the
fundamental property 
\begin{equation}
(D_{a}^{+}b)\cdot c=\left\{ 
\begin{array}{c}
c \\ 
a,b%
\end{array}%
\right\} ,\text{ for all }a,b,c\in \mathcal{V}(U).  \label{LCC.3a3}
\end{equation}%
Eq.(\ref{LCC.3a3}) follows immediately from Eq.(\ref{LCC.1}) once we change $%
c$ for $g^{-1}(c)$ and take into account the definitions given by Eq.(\ref%
{CHO.2}), Eq.(\ref{LCC.3}) and Eq.(\ref{LCC.3a1}).

We present now two remarkable properties of $\omega_{0}.\vspace{0.1in}$

\textbf{i.} For all $a,b,c\in\mathcal{V}(U)$ we have the \emph{cyclic
property} 
\begin{equation}
\omega_{0}(a)\underset{g}{\times}b\underset{g}{\cdot}c+\omega_{0}(b)\underset%
{g}{\times}c\underset{g}{\cdot}a+\omega_{0}(c)\underset{g}{\times }a\underset%
{g}{\cdot}b=0.  \label{LCC.3a}
\end{equation}

We show the cyclic property by recalling Eq.(\ref{LCC.N8}) used in the proof
of Eq.(\ref{LCC.1}). Indeed, we can write 
\begin{equation}
\omega_{0}(a)\underset{g}{\times}b\underset{g}{\cdot}c=\frac{1}{2}%
a\cdot((b\cdot\partial_{o}g)(c)-(c\cdot\partial_{o}g)(b)),  \label{LCC.M1}
\end{equation}
and, by cycling the letters $a,b$ and $c$, we get 
\begin{align}
\omega_{0}(b)\underset{g}{\times}c\underset{g}{\cdot}a & =\frac{1}{2}%
b\cdot((c\cdot\partial_{o}g)(a)-(a\cdot\partial_{o}g)(c)),  \label{LCC.M2} \\
\omega_{0}(c)\underset{g}{\times}a\underset{g}{\cdot}b & =\frac{1}{2}%
c\cdot((a\cdot\partial_{o}g)(b)-(b\cdot\partial_{o}g)(a)).  \label{LCC.M3}
\end{align}

Now, by adding Eqs.(\ref{LCC.M1}), (\ref{LCC.M2}) and (\ref{LCC.M3}) we get 
\begin{align*}
& \omega_{0}(a)\underset{g}{\times}b\underset{g}{\cdot}c+\omega _{0}(b)%
\underset{g}{\times}c\underset{g}{\cdot}a+\omega_{0}(c)\underset{g}{\times}a%
\underset{g}{\cdot}b \\
& =\frac{1}{2}(a\cdot(b\cdot\partial_{o}g)(c)-c\cdot(b\cdot\partial
_{o}g)(a))+\frac{1}{2}(b\cdot(c\cdot\partial_{o}g)(a)-a\cdot(c\cdot
\partial_{o}g)(b)) \\
& +\frac{1}{2}(c\cdot(a\cdot\partial_{o}g)(b)-b\cdot(a\cdot\partial
_{o}g)(c)).
\end{align*}
Then, by taking into account the symmetry property $(n\cdot\partial
_{o}g)^{\dagger}=n\cdot\partial_{o}g,$ the expected result immediately
follows.

\textbf{ii.} $\omega_{0}$ is just the $g$-\emph{gauge connection field}
associated to $(U,\lambda),$ i.e., 
\begin{equation}
\omega_{0}(a)=\frac{1}{2}\underset{g}{biv}[\lambda_{a}].  \label{LCC.3b}
\end{equation}

The proof of Eq.(\ref{LCC.3b}) uses the noticeable formulas: $(a\cdot
\partial _{o}\mathbf{t})\circ \mathbf{t}^{-1}+\mathbf{t}\circ (a\cdot
\partial _{o}\mathbf{t}^{-1})=0,$ for all \emph{non-singular} smooth $(1,1)$%
-extensor field $\mathbf{t}$, $(a\cdot \partial _{o}g^{-1})^{\dagger
}=a\cdot \partial _{o}g^{-1},$ $\partial _{n}\wedge (s(n))=0$, for all \emph{%
symmetric} $s\in ext_{1}^{1}(\mathcal{U}_{o})$, and $\partial _{n}\wedge
(B\times n)=-2B,$ where $B\in \bigwedge^{2}\mathcal{U}_{o}.$ A
straightforward calculation\footnote{%
Recall that $biv[t]=-\partial _{n}\wedge (t(n))$ and $\underset{g}{biv}%
[t]=biv[t\circ g^{-1}].$} allows us to get 
\begin{align*}
\underset{g}{biv}[\lambda _{a}]& =-\partial _{n}\wedge (\lambda _{a}\circ
g^{-1}(n)) \\
& =-\frac{1}{2}\partial _{n}\wedge g^{-1}\circ (a\cdot \partial _{o}g)\circ
g^{-1}(n)-\partial _{n}\wedge (\omega (a)\underset{g}{\times }g^{-1}(n)) \\
& =\frac{1}{2}\partial _{n}\wedge (a\cdot \partial _{o}g^{-1})(n)-\partial
_{n}\wedge (\omega (a)\times n), \\
& =0+2\omega _{0}(a).
\end{align*}

We present now three remarkable properties of the Levi-Civita connection
field $\lambda $ on $U$.\vspace{0.1in}

\textbf{i.} $\lambda$ is \emph{symmetric} with respect to the interchanging
of vector variables, i.e., 
\begin{equation}
\lambda(a,b)=\lambda(b,a).  \label{LCC.3c}
\end{equation}

To show this, let us take $a,b,c\in\mathcal{V}(U)$. Then, 
\begin{equation}
\lambda(a,b)\underset{g}{\cdot}c=\frac{1}{2}(a\cdot\partial_{o}g)(b)\cdot
c+\omega_{0}(a)\underset{g}{\times}b\underset{g}{\cdot}c,  \label{LCC.O1}
\end{equation}
and interchanging the letters $a$ and $b$ we have 
\begin{equation}
\lambda(b,a)\underset{g}{\cdot}c=\frac{1}{2}(b\cdot\partial_{o}g)(a)\cdot
c+\omega_{0}(b)\underset{g}{\times}a\underset{g}{\cdot}c.  \label{LCC.O2}
\end{equation}

Now, subtracting Eq.(\ref{LCC.O2}) from Eq.(\ref{LCC.O1}), and taking into
account Eq.(\ref{LCC.M3}) used in the proof of Eq.(\ref{LCC.3a}), we get 
\begin{equation*}
(\lambda(a,b)-\lambda(b,a))\underset{g}{\cdot}c=\omega_{0}(c)\underset{g}{%
\times}a\underset{g}{\cdot}b+\omega_{0}(a)\underset{g}{\times}b\underset{g}{%
\cdot}c-\omega_{0}(b)\underset{g}{\times}a\underset{g}{\cdot}c.
\end{equation*}
Then, by recalling the multivector identity $B\underset{g}{\times}v\underset{%
g}{\cdot}w=-B\underset{g}{\times}w\underset{g}{\cdot}v,$ where $%
B\in\bigwedge^{2}\mathcal{U}_{o}$ and $v,w\in\mathcal{U}_{o},$ and Eq.(\ref%
{LCC.3a}), the required result immediately follows by the non-degeneracy of
the $g$-scalar product.

\textbf{ii.} The $g$-\emph{symmetric} and $g$-\emph{skew symmetric} parts of 
$\lambda_{a},$ namely $\lambda_{a\pm(g)}=\dfrac{1}{2}(\lambda_{a}\pm
\lambda_{a}^{\dagger(g)}),$ are given by 
\begin{align}
\lambda_{a+(g)}(b) & =\dfrac{1}{2}g^{-1}\circ(a\cdot\partial_{o}g)(b),
\label{LCC.3d} \\
\lambda_{a-(g)}(b) & =\omega_{0}(a)\underset{g}{\times}b.  \label{LCC.3f}
\end{align}

To show these important formulas we calculate first the metric adjoint of $%
\lambda_{a},$ namely $\lambda_{a}^{\dagger(g)}$, by using the fundamental
property \cite{2} of the metric adjoint operator $\left. {}\right.
^{\dagger(g)}$. By recalling the symmetry property $(a\cdot\partial
_{o}g)^{\dagger}=a\cdot\partial_{o}g$ and the multivector identity $B%
\underset{g}{\times}v\underset{g}{\cdot}w=-B\underset{g}{\times}w\underset{g}%
{\cdot}v,$ where $B\in\bigwedge^{2}\mathcal{U}_{o}$ and $v,w\in\mathcal{U}%
_{o},$ we can write that 
\begin{align*}
\lambda_{a}^{\dagger(g)}(b)\underset{g}{\cdot}c & =b\underset{g}{\cdot }%
\lambda_{a}(c) \\
& =\frac{1}{2}b\cdot(a\cdot\partial_{o}g)(c)-\omega_{0}(a)\underset{g}{\times%
}b\underset{g}{\cdot}c \\
& =(\frac{1}{2}g^{-1}\circ(a\cdot\partial_{o}g)(b)-\omega_{0}(a)\underset{g}{%
\times}b)\underset{g}{\cdot}c,
\end{align*}
hence, by the non-degeneracy of the $g$-scalar product, it follows that 
\begin{equation*}
\lambda_{a}^{\dagger(g)}(b)=\frac{1}{2}g^{-1}\circ(a\cdot\partial
_{o}g)(b)-\omega_{0}(a)\underset{g}{\times}b.
\end{equation*}
Now, we can get that 
\begin{align*}
\lambda_{a}(b)+\lambda_{a}^{\dagger(g)}(b) & =g^{-1}\circ(a\cdot\partial
_{o}g)(b), \\
\lambda_{a}(b)-\lambda_{a}^{\dagger(g)}(b) & =2\omega_{0}(a)\underset{g}{%
\times}b.
\end{align*}

\textbf{iii.} The \emph{generalized }of $\lambda _{a},$ namely $\Lambda
_{a}, $ is given by the following formula 
\begin{equation}
\Lambda _{a}(X)=\dfrac{1}{2}\underline{g}^{-1}\circ (a\cdot \partial _{o}%
\underline{g})(X)+\omega _{0}(a)\underset{g}{\times }X,  \label{LCC.4}
\end{equation}%
where $\underline{g}$ and $\underline{g}^{-1}$ are \emph{extended} of $g$
and extended of $g^{-1}$, respectively.

The above property is an immediate consequence of using the noticeable
formulas: $\mathbf{t}^{-1}\circ (a\cdot \partial _{o}\mathbf{t})(\partial
_{n})\wedge (n\lrcorner X)=\underline{\mathbf{t}}^{-1}\circ (a\cdot \partial
_{o}\underline{\mathbf{t}})(X),$ for all \emph{non-singular} smooth $(1,1)$%
-extensor field $\mathbf{t},$ and $(B\underset{g}{\times }\partial
_{n})\wedge (n\lrcorner X)=B\underset{g}{\times }X,$ where $B\in
\bigwedge^{2}\mathcal{U}_{o}$ and $X\in \bigwedge \mathcal{U}_{o}$.

From the above, we suspect (and this is indeed the case, see next section)
that a a particular $\lambda _{a}$ which is metric compatible and symmetric
may represent in $U$ the effects of the restriction in $U$ of the
Levi-Civita connection of $%
\slg%
$ on $M$.

\section{ Covariant Derivatives and Metric Compatibility}

\subsection{Geometric Structure}

A \emph{parallelism structure }$(U,\gamma)$ is said to be \emph{compatible}
with a \emph{metric structure} $(U,g)$ if and only if 
\begin{equation}
\gamma_{a+(g)}=\frac{1}{2}g^{-1}\circ(a\cdot\partial_{o}g),  \label{GS.1}
\end{equation}
i.e., $g\circ\gamma_{a}+\gamma_{a}^{\dagger}\circ g=a\cdot\partial_{o}g$.

Sometimes for abuse of language we will say that $\gamma $ is \emph{metric
compatible} (or $g$-\emph{compatible}, for short).

\subsubsection{Geometrical Structure in $U$}

The open set $U$ equipped with a connection field $\gamma $ and a metric
extensor field $g,$ namely $(U,\gamma ,g),$ such that $(U,\gamma )$ is
compatible with $(U,g),$ will be said to be a \emph{geometric structure} on $%
U.$

The Levi-Civita parallelism structure $(U,\lambda ),$ according to Eq.(\ref%
{LCC.3d}), is compatible with the metric structure $(U,g),$ i.e., $\lambda $
is $g$-compatible. It immediately follows that $(U,\lambda ,g)$ is a
well-defined geometric structure on $U.$

Using the above results we have in any geometric structure $(U,\gamma ,g)$
an important result.

\textbf{Theorem 1. }There exists a smooth $(1,2)$-extensor field on $U,$
namely $\omega,$ such that 
\begin{equation}
\gamma_{a}(b)=\frac{1}{2}g^{-1}\circ(a\cdot\partial_{o}g)(b)+\omega (a)%
\underset{g}{\times}b.  \label{GS.1a}
\end{equation}

We give now three properties involving $\gamma_{a}$ and $\omega.\vspace {%
0.1in}$

\textbf{i.} The $g$-\emph{symmetric} and $g$-\emph{skew-symmetric parts} of $%
\gamma_{a},$ namely $\gamma_{a\pm(g)}=\dfrac{1}{2}(\gamma_{a}\pm\gamma
_{a}^{\dagger}),$ are given by 
\begin{align}
\gamma_{a+(g)}(b) & =\frac{1}{2}g^{-1}\circ(a\cdot\partial_{o}g)(b),
\label{GS.1b} \\
\gamma_{a-(g)}(b) & =\omega(a)\underset{g}{\times}b.  \label{GS.1c}
\end{align}

\textbf{ii.} The \emph{generalized} of $\gamma_{a},$ namely $\Gamma_{a},$ is
given by 
\begin{equation}
\Gamma_{a}(X)=\frac{1}{2}\underline{g}^{-1}\circ(a\cdot\partial_{o}%
\underline{g})(X)+\omega(a)\underset{g}{\times}X.  \label{GS.1d}
\end{equation}

\textbf{iii.} $\omega$ is just the $g$-\emph{gauge connection field}
associated to $(U,\gamma),$ i.e., 
\begin{equation}
\omega(a)=\frac{1}{2}\underset{g}{biv}[\gamma_{a}].  \label{GS.1e}
\end{equation}

\subsection{Metric Compatibility}

The pair of $a$-\emph{DCDO's }associated to $(U,\gamma),$ namely $(\nabla
_{a}^{+},\nabla_{a}^{-}),$ is said to be \emph{metric} \emph{compatible} (or 
$g$-\emph{compatible}, for short) if and only if 
\begin{equation}
\nabla_{a}^{++}g=0,  \label{MCD.1}
\end{equation}
or equivalently, 
\begin{equation}
\nabla_{a}^{--}g^{-1}=0.  \label{MCD.1a}
\end{equation}

We emphasize that Eq.(\ref{MCD.1}) and Eq.(\ref{MCD.1a}) are completely
equivalent to each other. It follows from the remarkable formula $(\nabla
_{a}^{++}\mathbf{t})\circ \mathbf{t}^{-1}+\mathbf{t}\circ (\nabla _{a}^{--}%
\mathbf{t})=0$, valid for all non-singular smooth $(1,1)$-extensor field $%
\tau .$

$(U,\gamma)$ is compatible with $(U,g)$ if and only if $(\nabla_{a}^{+},%
\nabla_{a}^{-})$ is metric compatible.

Indeed, let us take $b\in\mathcal{V}(U),$ we have that 
\begin{align}
(\nabla_{a}^{++}g)(b) & =\nabla_{a}^{-}g(b)-g(\nabla_{a}^{+}b)  \notag \\
& =a\cdot\partial_{o}g(b)-\gamma_{a}^{\dagger}\circ g(b)-g(a\cdot\partial
_{o}b)-g\circ\gamma_{a}(b),  \notag \\
& =(a\cdot\partial_{o}g)(b)-g\circ\gamma_{a}(b)-\gamma_{a}^{\dagger}\circ
g(b).  \label{MCD.1b}
\end{align}
Now, if $\gamma$ is $g$-compatible, by using Eq.(\ref{GS.1}) into Eq.(\ref%
{MCD.1b}), it follows that $(\nabla_{a}^{++}g)(b)=0,$ i.e., $%
(\nabla_{a}^{+},\nabla_{a}^{-})$ is $g$-compatible. And, if $%
(\nabla_{a}^{+},\nabla_{a}^{-})$ is $g$-compatible, by using Eq.(\ref{MCD.1}%
) into Eq.(\ref{MCD.1b}), we get that $g\circ\gamma_{a}(b)+\gamma_{a}^{%
\dagger}\circ g(b)=(a\cdot\partial_{o}g)(b),$ i.e., $\gamma$ is $g$%
-compatible.

We now present some basic properties which are satisfied by a $g$-compatible
pair of $a$-\emph{DCDO's,} namely $(\mathcal{D}_{a}^{+},\mathcal{D}_{a}^{-})$%
.\vspace{0.1in}

\textbf{i.} For any $(\mathcal{D}_{a}^{+},\mathcal{D}_{a}^{-})$ we have 
\begin{align}
\mathcal{D}_{a}^{++}\underline{g} & =0,  \label{MCD.2} \\
\mathcal{D}_{a}^{\_\_}\underline{g}^{-1} & =0,  \label{MCD.2a}
\end{align}
where $\underline{g}$ and $\underline{g}^{-1}$are the so-called \emph{%
extended} of $g$ and $g^{-1},$ respectively.

In order to prove the first statement we only need to check that for all $%
f\in\mathcal{S}(U)$ and $b_{1},\ldots,b_{k}\in\mathcal{V}(U)$ 
\begin{equation*}
(\mathcal{D}_{a}^{++}\underline{g})(f)=0\text{ and }(\mathcal{D}_{a}^{++}%
\underline{g})(b_{1}\wedge\ldots\wedge b_{k})=0.
\end{equation*}

But, by using the fundamental property $\underline{g}(f)=f,$ we get 
\begin{equation*}
(\mathcal{D}_{a}^{++}\underline{g})(f)=\mathcal{D}_{a}^{-}\underline {g}(f)-%
\underline{g}(\mathcal{D}_{a}^{+}f)=\mathcal{D}_{a}^{-}f-\underline {g}%
(a\cdot\partial_{o}f)=a\cdot\partial_{o}f-a\cdot\partial_{o}f=0.
\end{equation*}

And, by using the fundamental property $\underline{g}(b_{1}\wedge\ldots%
\wedge b_{k})=g(b_{1})\wedge\ldots\wedge g(b_{k})$ we get 
\begin{align*}
& (\mathcal{D}_{a}^{++}\underline{g})(b_{1}\wedge\ldots\wedge b_{k}) \\
& =\mathcal{D}_{a}^{-}\underline{g}(b_{1}\wedge\ldots\wedge b_{k})-%
\underline{g}(\mathcal{D}_{a}^{+}(b_{1}\wedge\ldots\wedge b_{k})) \\
& =\mathcal{D}_{a}^{-}g(b_{1})\wedge\ldots
g(b_{k})+\cdots+g(b_{1})\wedge\ldots\mathcal{D}_{a}^{-}g(b_{k}) \\
& -g(\mathcal{D}_{a}^{+}b_{1})\wedge\ldots
g(b_{k})-\ldots-g(b_{1})\wedge\ldots g(\mathcal{D}_{a}^{+}b_{k}) \\
& =(\mathcal{D}_{a}^{++}g)(b_{1})\wedge\ldots
g(b_{k})+\cdots+g(b_{1})\wedge\ldots(\mathcal{D}_{a}^{++}g)(b_{k})=0,
\end{align*}
since $\mathcal{D}_{a}^{++}g=0.$

The second statement can be proved analogously.

\textbf{ii. }$(\mathcal{D}_{a}^{+},\mathcal{D}_{a}^{-})$ satisfies the
fundamental property 
\begin{equation}
\mathcal{D}_{a}^{-}\underline{g}(X)=\underline{g}(\mathcal{D}_{a}^{+}X).
\label{MCD.3}
\end{equation}

Indeed, by Eq.(\ref{MCD.2}) we have that for all $X\in\mathcal{M}(U)$ 
\begin{align*}
(\mathcal{D}_{a}^{++}\underline{g})(X) & =0, \\
\mathcal{D}_{a}^{-}\underline{g}(X)-\underline{g}(\mathcal{D}_{a}^{+}X) & =0.
\end{align*}

\textbf{iii. }\emph{Ricci-like theorems}. Let $X,Y$ be a smooth multivector
fields. Then, 
\begin{align}
a\cdot\partial_{o}(X\underset{g}{\cdot}Y) & =(\mathcal{D}_{a}^{+}X)\underset{%
g}{\cdot}Y+X\underset{g}{\cdot}(\mathcal{D}_{a}^{+}Y),  \label{MCD.4} \\
a\cdot\partial_{o}(X\underset{g^{-1}}{\cdot}Y) & =(\mathcal{D}_{a}^{-}X)%
\underset{g^{-1}}{\cdot}Y+X\underset{g^{-1}}{\cdot}(\mathcal{D}_{a}^{-}Y).
\label{MCD.4b}
\end{align}

The proof is as follows. As we know, $(\mathcal{D}_{a}^{+},\mathcal{D}%
_{a}^{-})$ must satisfy the fundamental property 
\begin{equation*}
(\mathcal{D}_{a}^{+}X)\cdot Y+X\cdot(\mathcal{D}_{a}^{-}Y)=a\cdot\partial
_{o}(X\cdot Y).
\end{equation*}
Then, by substituting $Y$ for $\underline{g}(Y)$ and using Eq.(\ref{MCD.3}),
the first statement follows immediately. To prove get the second statement
it is enough to substitute $Y$ for $\underline{g}^{-1}(Y)$ and once again
use Eq.(\ref{MCD.3}).

\textbf{iii. }$(\mathcal{D}_{a}^{+},\mathcal{D}_{a}^{-})$ satisfies \emph{%
Leibnitz-like rules} for all of the $g$ and $g^{-1}$ suitable products%
\footnote{%
Here, as in \cite{2} $\underset{g}{\ast}$ means any product, either $%
(\wedge),$ $(\underset{g}{\cdot}),$ $(\underset{g}{\lrcorner },\underset{g}{%
\llcorner})$ or $(g$-\emph{Clifford product}$).$ Analogously for $\underset{%
g^{-1}}{\ast}.$}, namely $\underset{g}{\ast}$ and $\underset{g^{-1}}{\ast},$
of smooth multivector fields 
\begin{align}
\mathcal{D}_{a}^{+}(X\underset{g}{\ast}Y) & =(\mathcal{D}_{a}^{+}X)\underset{%
g}{\ast}Y+X\underset{g}{\ast}(\mathcal{D}_{a}^{+}Y),  \label{MCD.5} \\
\mathcal{D}_{a}^{-}(X\underset{g^{-1}}{\ast}Y) & =(\mathcal{D}_{a}^{-}X)%
\underset{g^{-1}}{\ast}Y+X\underset{g^{-1}}{\ast}(\mathcal{D}_{a}^{-}Y).
\label{MCD.5a}
\end{align}

We prove only the first statement. The other proof is similar.

Firstly, if $\underset{g}{\ast }$ is just$(\wedge ),$ then Eq.(\ref{MCD.5})
is nothing more than the Leibniz rule for the exterior product of smooth
multivector fields i.e., 
\begin{equation}
\mathcal{D}_{a}^{+}(X\wedge Y)=(\mathcal{D}_{a}^{+}X)\wedge Y+X\wedge (%
\mathcal{D}_{a}^{+}Y).  \label{MCD.5i}
\end{equation}%
As we know it is true.

Secondly, if $\underset{g}{\ast}$ is $(\underset{g}{\cdot})$, then $\mathcal{%
D}_{a}^{+}(X\underset{g}{\ast}Y)=a\cdot\partial_{o}(X\underset{g}{\cdot}Y)$
and it follows that Eq.(\ref{MCD.5}) is nothing more than the Ricci-like
theorem for $\mathcal{D}_{a}^{+}.$

In order to prove Eq.(\ref{MCD.5}) whenever $\underset{g}{\ast}$ is either $(%
\underset{g}{\lrcorner})$ or $(\underset{g}{\llcorner}),$ we use the
identities $(X\underset{g}{\lrcorner}Y)\underset{g}{\cdot}Z=Y\underset{g}{%
\cdot}(\widetilde{X}\wedge Z)$ and $(X\underset{g}{\llcorner}Y)\underset{g}{%
\cdot}Z=X\underset{g}{\cdot}(Z\wedge\widetilde{Y}),$ for all $X,Y,Z\in%
\mathcal{M}(U),$ and Eq.(\ref{MCD.4}) and Eq.(\ref{MCD.5i}). For instance,
for the left $g$-contracted product we can indeed write%
\begin{equation*}
a\cdot\partial_{o}((X\underset{g}{\lrcorner}Y)\underset{g}{\cdot}%
Z)=a\cdot\partial_{o}(Y\underset{g}{\cdot}(\widetilde{X}\wedge Z))
\end{equation*}
and then 
\begin{align*}
& (\mathcal{D}_{a}^{+}(X\underset{g}{\lrcorner}Y))\underset{g}{\cdot }Z+(X%
\underset{g}{\lrcorner}Y)\underset{g}{\cdot}(\mathcal{D}_{a}^{+}Z) \\
& =(\mathcal{D}_{a}^{+}Y)\underset{g}{\cdot}(\widetilde{X}\wedge Z)+Y%
\underset{g}{\cdot}((\mathcal{D}_{a}^{+}\widetilde{X})\wedge Z)+Y\underset{g}%
{\cdot}(\widetilde{X}\wedge(\mathcal{D}_{a}^{+}Z))
\end{align*}
Also,%
\begin{equation*}
(\mathcal{D}_{a}^{+}(X\underset{g}{\lrcorner}Y))\underset{g}{\cdot }Z=((%
\mathcal{D}_{a}^{+}X)\underset{g}{\lrcorner}Y+X\underset{g}{\lrcorner }(%
\mathcal{D}_{a}^{+}Y))\underset{g}{\cdot}Z.
\end{equation*}

Hence, by the non-degeneracy of the $g$-scalar product, the Leibniz rule for 
$(\underset{g}{\lrcorner})$ immediately follows, i.e., 
\begin{equation}
\mathcal{D}_{a}^{+}(X\underset{g}{\lrcorner}Y)=(\mathcal{D}_{a}^{+}X)%
\underset{g}{\lrcorner}Y+X\underset{g}{\lrcorner}(\mathcal{D}_{a}^{+}Y).
\label{MCD.5ii}
\end{equation}

In order to prove Eq.(\ref{MCD.5}) whenever $\underset{g}{\ast}$ means $(%
\underset{g}{}),$ we only need to check that for all $f\in\mathcal{S}(U)$
and $b_{1},\ldots,b_{k}\in\mathcal{V}(U)$ 
\begin{align}
\mathcal{D}_{a}^{+}(f\underset{g}{}X) & =(\mathcal{D}_{a}^{+}f)\underset{g}{}%
X+f\underset{g}{}(\mathcal{D}_{a}^{+}X),  \label{MCD.5iii} \\
\mathcal{D}_{a}^{+}(b_{1}\underset{g}{}\ldots\underset{g}{}b_{k}\underset{g}{%
}X) & =(\mathcal{D}_{a}^{+}(b_{1}\underset{g}{}\ldots\underset{g}{}b_{k}))%
\underset{g}{}X  \notag \\
& +b_{1}\underset{g}{}\ldots\underset{g}{}b_{k}\underset{g}{}(\mathcal{D}%
_{a}^{+}X).  \label{MCD.5iv}
\end{align}

The verification of Eq.(\ref{MCD.5iii}) is trivial.

To verify Eq.(\ref{MCD.5iv}) we will use \emph{complete induction }over the $%
k$ smooth vector fields $b_{1},\ldots,b_{k}.$

Let us take $b\in\mathcal{V}(U),$ by using Eq.(\ref{MCD.5ii}) and Eq.(\ref%
{MCD.5i}), we have 
\begin{align}
\mathcal{D}_{a}^{+}(b\underset{g}{}X) & =\mathcal{D}_{a}^{+}(b\underset{g}{%
\lrcorner}X)+\mathcal{D}_{a}^{+}(b\wedge X)  \notag \\
& =(\mathcal{D}_{a}^{+}b)\underset{g}{\lrcorner}X+b\underset{g}{\lrcorner }(%
\mathcal{D}_{a}^{+}X)+(\mathcal{D}_{a}^{+}b)\wedge X+b\wedge(\mathcal{D}%
_{a}^{+}X),  \notag \\
\mathcal{D}_{a}^{+}(b\underset{g}{}X) & =(\mathcal{D}_{a}^{+}b)\underset{g}{}%
X+b\underset{g}{}(\mathcal{D}_{a}^{+}X).  \label{MCD.5v}
\end{align}

Now, let us $b_{1},\ldots,b_{k},b_{k+1}\in\mathcal{V}(U).$ By using twice
the inductive hypothesis and Eq.(\ref{MCD.5v}) we can write 
\begin{align*}
& \mathcal{D}_{a}^{+}(b_{1}\underset{g}{}\ldots\underset{g}{}b_{k}\underset{g%
}{}b_{k+1}\underset{g}{}X) \\
& =(\mathcal{D}_{a}^{+}(b_{1}\underset{g}{}\ldots\underset{g}{}b_{k}))%
\underset{g}{}b_{k+1}\underset{g}{}X+b_{1}\underset{g}{}\ldots \underset{g}{}%
b_{k}\underset{g}{}(\mathcal{D}_{a}^{+}(b_{k+1}\underset{g}{}X)) \\
& =(\mathcal{D}_{a}^{+}(b_{1}\underset{g}{}\ldots\underset{g}{}b_{k}))%
\underset{g}{}b_{k+1}\underset{g}{}X+b_{1}\underset{g}{}\ldots \underset{g}{}%
b_{k}\underset{g}{}(\mathcal{D}_{a}^{+}b_{k+1})\underset{g}{}X \\
& +b_{1}\underset{g}{}\ldots\underset{g}{}b_{k}\underset{g}{}b_{k+1}\underset%
{g}{}(\mathcal{D}_{a}^{+}X) \\
& =((\mathcal{D}_{a}^{+}(b_{1}\underset{g}{}\ldots\underset{g}{}b_{k}))%
\underset{g}{}b_{k+1}+b_{1}\underset{g}{}\ldots\underset{g}{}b_{k}\underset{g%
}{}(\mathcal{D}_{a}^{+}b_{k+1}))\underset{g}{}X \\
& +b_{1}\underset{g}{}\ldots\underset{g}{}b_{k}\underset{g}{}b_{k+1}\underset%
{g}{}(\mathcal{D}_{a}^{+}X) \\
& =(\mathcal{D}_{a}^{+}(b_{1}\underset{g}{}\ldots\underset{g}{}b_{k}\underset%
{g}{}b_{k+1}))\underset{g}{}X+b_{1}\underset{g}{}\ldots\underset{g}{}b_{k}%
\underset{g}{}b_{k+1}\underset{g}{}(\mathcal{D}_{a}^{+}X).
\end{align*}

From the theory of extensors developed in (\cite{2}) it immediately follows
that given any metric $g$, there is a \emph{non-singular} smooth $(1,1)$%
-extensor field $h$ such that 
\begin{equation}
g=h^{\dagger}\circ\eta\circ h,  \label{MCD.6}
\end{equation}
where $\eta$ is an \emph{pseudo-orthogonal metric field} with the same
signature as $g$. Such $h$ will be said to be a \emph{gauge metric field}
for $g$.

We now ask if there is any noticeable relationship between a $g$-compatible
pair of $a$-\emph{DCDO's} and a $\eta $-compatible pair of $a$-\emph{DCDO's.}
The answer is YES.

\textbf{Theorem 2. }Let $h$ be a gauge metric field for $g$. For any $g$%
-compatible pair of $a$-\emph{DCDO's,} namely $(_{g}\mathcal{D}_{a}^{+}$ $%
,_{g}\mathcal{D}_{a}^{-}),$ there exists an unique $\eta$-compatible pair of 
$a$-\emph{DCDO's,} namely $(_{\eta}\mathcal{D}_{a}^{+},_{\eta}\mathcal{D}%
_{a}^{-}),$ such that\footnote{%
Recall that $h^{*}=(h^{-1})^{\dagger }=(h^{\dagger})^{-1},$ and $\underline{h%
}^{-1}=(\underline{h})^{-1}=\underline{(h^{-1})}$ and $\underline{h}%
^{\dagger}=(\underline{h})^{\dagger }=\underline{(h^{\dagger})}.$} 
\begin{align}
\underline{h}(_{g}\mathcal{D}_{a}^{+}X) & =\text{ }_{\eta}\mathcal{D}_{a}^{+}%
\underline{h}(X),  \label{MCD.7a} \\
\underline{h}^{\ast}(_{g}\mathcal{D}_{a}^{-}X) & =\text{ }_{\eta}\mathcal{D}%
_{a}^{-}\underline{h}^{\ast}(X).  \label{MCD.7b}
\end{align}
And reciprocally, given any $\eta$-compatible pair of $a$-\emph{DCDO's}, say 
$(_{\eta}\mathcal{D}_{a}^{+},_{\eta}\mathcal{D}_{a}^{-}),$ there is an
unique $g$-compatible pair of $a$-\emph{DCDO's,} say $(_{g}\mathcal{D}%
_{a}^{+},_{g}\mathcal{D}_{a}^{-}),$ such that the above formulas are
satisfied.

\textbf{Proof}

Given $(_{g}\mathcal{D}_{a}^{+},_{g}\mathcal{D}_{a}^{-}),$ since $h$ is a
non-singular smooth $(1,1)$-extensor field, we can indeed construct a
well-defined pair of $a$-\emph{DCDO's,} namely $(_{h}\mathcal{D}_{a}^{+},_{h}%
\mathcal{D}_{a}^{-}),$ by the following formulas 
\begin{equation*}
_{h}\mathcal{D}_{a}^{+}X=\underline{h}(_{g}\mathcal{D}_{a}^{+}\underline {h}%
^{-1}(X))\text{ and }_{h}\mathcal{D}_{a}^{-}X=\underline{h}^{\ast}(_{h}%
\mathcal{D}_{a}^{+}\underline{h}^{\dagger}(X)).
\end{equation*}
As defined above, $(_{h}\mathcal{D}_{a}^{+},_{h}\mathcal{D}_{a}^{-})$ is the 
$h$-\emph{deformation} of $(_{g}\mathcal{D}_{a}^{+},_{g}\mathcal{D}_{a}^{-})$%
.

But, it is obvious that $_{h}\mathcal{D}_{a}^{+}$ and $_{h}\mathcal{D}%
_{a}^{-}$ as defined above satisfy in fact Eq.(\ref{MCD.7a}) and Eq.(\ref%
{MCD.7b}), i.e., $\underline{h}(_{g}\mathcal{D}_{a}^{+}X)=$ $_{h}\mathcal{D}%
_{a}^{+}\underline{h}(X)$ and $\underline{h}^{\ast}(_{g}\mathcal{D}%
_{a}^{-}X)=$ $_{\eta}\mathcal{D}_{a}^{-}\underline{h}^{\ast}(X).$

In order to check the $\eta$-compatibility of $(_{h}\mathcal{D}_{a}^{+},_{h}%
\mathcal{D}_{a}^{-}),$ we can write 
\begin{align*}
(_{h}\mathcal{D}_{a}^{++}\eta)(b) & =\text{ }_{h}\mathcal{D}%
_{a}^{-}\eta(b)-\eta(_{h}\mathcal{D}_{a}^{+}b) \\
& =h^{\ast}(_{g}\mathcal{D}_{a}^{-}h^{\dagger}\circ\eta(b))-\eta\circ h(_{g}%
\mathcal{D}_{a}^{+}h^{-1}(b)) \\
& =h^{\ast}(_{g}\mathcal{D}_{a}^{-}h^{\dagger}\circ\eta\circ h\circ
h^{-1}(b)-h^{\dagger}\circ\eta\circ h(_{g}\mathcal{D}_{a}^{+}h^{-1}(b))) \\
& =h^{\ast}(_{g}\mathcal{D}_{a}^{-}g(h^{-1}(b))-g(_{g}\mathcal{D}%
_{a}^{+}h^{-1}(b))), \\
& =h^{\ast}(_{g}\mathcal{D}_{a}^{++}g)(h^{-1}(b)).
\end{align*}
This implies that $_{h}\mathcal{D}_{a}^{++}\eta=h^{\ast}\circ(_{g}\mathcal{D}%
_{a}^{++}g)\circ h^{-1}$. Then, since $_{g}\mathcal{D}_{a}^{++}g=0,$ it
follows that $_{h}\mathcal{D}_{a}^{++}\eta=0,$ i.e., $(_{h}\mathcal{D}%
_{a}^{+},_{h}\mathcal{D}_{a}^{-})$ is $\eta$-compatible.

Now, if there exists another $\eta$-compatible pair $(_{\eta}\mathcal{D}%
_{a}^{\prime+},_{\eta}\mathcal{D}_{a}^{\prime-})$, which satisfies Eq.(\ref%
{MCD.7a}) and Eq.(\ref{MCD.7b}), i.e., 
\begin{equation*}
\underline{h}(_{g}\mathcal{D}_{a}^{+}X)=\text{ }_{\eta}\mathcal{D}%
_{a}^{\prime+}\underline{h}(X)\text{ and }\underline{h}^{\ast}(_{g}\mathcal{D%
}_{a}^{-}X)=\text{ }_{\eta}\mathcal{D}_{a}^{\prime-}\underline{h}^{\ast}(X).
\end{equation*}
Then, by substituting $X$ for $\underline{h}^{-1}(X)$ in the first one and $%
X $ for $\underline{h}^{\dagger}(X)$ in the second one, it follows that $%
_{\eta }\mathcal{D}_{a}^{\prime+}=$ $_{\eta}\mathcal{D}_{a}^{+}$ and $_{\eta
}\mathcal{D}_{a}^{\prime-}=$ $_{\eta}\mathcal{D}_{a}^{-}$.

So the existence and uniqueness are proved. Such a $\eta$-compatible pair of 
$a$-\emph{DCDO's} satisfying Eq.(\ref{MCD.7a}) and Eq.(\ref{MCD.7b}) is just
the $h$-deformation of the $g$-compatible pair of $a$-\emph{DCDO's.}

By following analogous steps we can also prove that such a $g$-compatible
pair of $a$-\emph{DCDO's} satisfying Eq.(\ref{MCD.7a}) and Eq.(\ref{MCD.7b})
is just the $h^{-1}$-deformation of the $\eta$-compatible pair of $a$-\emph{%
DCDO's}.$\blacksquare$

In the fourth paper of this series we study the relation between the
curvature and torsion tensors of a pair $(_{\eta }\mathcal{D}_{a}^{+},_{\eta
}\mathcal{D}_{a}^{-})$ and its deformation $(_{g}\mathcal{D}_{a}^{+},_{g}%
\mathcal{D}_{a}^{-})$.

\section{Derivative Operators in Metric and Geometric Structures}

\subsection{Ordinary Hodge Coderivatives}

Let $U$ be an open subset of $U_{o},$ and let $(U,g)$ be a \emph{metric
structure} on $U.$ Let us take any pair of \emph{reciprocal frame fields} on 
$U,$ say $(\{e_{\mu}\},\{e^{\mu}\}),$ i.e., $e_{\mu}\cdot
e^{\nu}=\delta_{\mu }^{\nu}.$ In particular, the \emph{fiducial frame field}%
, namely $\{b_{\mu }\},$ due to its \emph{orthonormality}, i.e., $%
b_{\mu}\cdot b_{\nu}=\delta_{\mu\nu},$ is a \emph{self-reciprocal frame field%
}, i.e., $b^{\mu }=b_{\mu}.$ We note also that $\{b_{\mu}\}$ is an \emph{%
ordinarily constant frame field} on $U,$ i.e., 
\begin{equation}
a\cdot\partial_{o}b_{\mu}=0,\text{ for each }\mu=1,\ldots,n.  \label{OHD.1}
\end{equation}

Associated to $(\{e_{\mu}\},\{e^{\mu}\}),$ the smooth pseudoscalar field on $%
U,$ namely $\tau,$ defined by 
\begin{equation}
\tau=\sqrt{e_{\wedge}\cdot e_{\wedge}}e^{\wedge},  \label{OHD.2}
\end{equation}
where $e_{\wedge}=e_{1}\wedge\ldots\wedge e_{n}\in\mathcal{M}^{n}(U)$ and $%
e^{\wedge}=e^{1}\wedge\ldots\wedge e^{n}\in\mathcal{M}^{n}(U),$ is said to
be the \emph{standard volume pseudoscalar field }for the \emph{local
coordinate system} $(U_{o},\phi_{o}).$

Such $\tau\in\mathcal{M}^{n}(U)$ has the fundamental property 
\begin{equation}
\tau\cdot\tau=\tau\lrcorner\widetilde{\tau}=\tau\widetilde{\tau}=1.
\label{OHD.2a}
\end{equation}
It follows from the obvious result $e_{\wedge}\cdot e^{\wedge}=1.$

From Eq.(\ref{OHD.2a}) we can get an expansion formula for smooth
pseudoscalar fields on $U,$ i.e., 
\begin{equation}
I=(I\cdot\tau)\tau.  \label{OHD.2b}
\end{equation}

In particular, because of the obvious properties $b_{\wedge}\cdot b_{\wedge
}=1$ and $b^{\wedge}=b_{\wedge},$ the standard volume pseudoscalar field
associated to $\{b_{\mu}\}$ is just $b_{\wedge}.$ It will be called the 
\emph{canonical volume pseudoscalar field} for $(U_{o},\phi_{o}).$

We emphasize that $b_{\wedge }$ is an \emph{ordinarily constant }smooth
pseudoscalar field on $U,$ i.e., 
\begin{equation}
a\cdot \partial _{o}b_{\wedge }=0.  \label{OHD.3}
\end{equation}%
Eq.(\ref{OHD.3}) can be proved by using Eq.(\ref{OHD.1}) and the Leibniz
rule for the exterior product of smooth multivector fields.

By using Eq.(\ref{OHD.3}) and the general Leibniz rule for $a\cdot \partial
_{o}$ we can deduce the remarkable property 
\begin{equation}
a\cdot \partial _{o}(b_{\wedge }\ast X)=b_{\wedge }\ast (a\cdot \partial
_{o}X),  \label{OHD.4}
\end{equation}%
where $\ast $ means either exterior product or any canonical product of
smooth multivector fields.

On the other hand we have that all $(\{e_{\mu}\},\{e^{\mu}\})$ must be
necessarily an \emph{extensor-deformation} of $\{b_{\mu}\}$. This statement
means that there exists a non-singular smooth $(1,1)$-extensor field on $U,$
say $\varepsilon,$ such that 
\begin{align}
e_{\mu} & =\varepsilon(b_{\mu}),  \label{OHD.5a} \\
e^{\mu} & =\varepsilon^{*}(b_{\mu}),\text{ for each }\mu=1,\ldots,n.
\label{OHD.5b}
\end{align}

Then, by putting Eq.(\ref{OHD.5a}) and Eq.(\ref{OHD.5b}) into Eq.(\ref{OHD.2}%
), we have that 
\begin{equation*}
\tau=(\underline{\varepsilon}(b_{\wedge})\cdot\underline{\varepsilon }%
(b_{\wedge}))^{1/2}\underline{\varepsilon}^{\ast}(b_{\wedge})=(\left.
\det\right. ^{2}[\varepsilon]b_{\wedge}\cdot b_{\wedge})^{1/2}\left.
\det\right. ^{-1}[\varepsilon]b_{\wedge}=sgn(\det[\varepsilon])b_{\wedge},
\end{equation*}
i.e., 
\begin{equation}
\tau=\pm b_{\wedge}.  \label{OHD.6}
\end{equation}

From Eq.(\ref{OHD.3}) and Eq.(\ref{OHD.4}), by taking into account Eq.(\ref%
{OHD.6}), we have two remarkable properties 
\begin{align}
a\cdot\partial_{o}\tau & =0,  \label{OHD.6a} \\
a\cdot\partial_{o}(\tau\ast X) & =\tau\ast(a\cdot\partial_{o}X),\text{ for
all }X\in\mathcal{M}(U).  \label{OHD.6b}
\end{align}

Associated to $(\{e_{\mu}\},\{e^{\mu}\}),$ the smooth pseudoscalar field on $%
U,$ namely $\underset{g}{\tau},$ defined by 
\begin{equation}
\underset{g}{\tau}=\sqrt{\left| e_{\wedge}\underset{g}{\cdot}e_{\wedge
}\right| }e^{\wedge}=\sqrt{\left| \det[g]\right| }\tau,  \label{OHD.7}
\end{equation}
will be said to be a \emph{metric volume pseudoscalar field} for $%
(U_{o},\phi_{o}).$

Such $\underset{g}{\tau}$ $\in\mathcal{M}^{n}(U)$ satisfies the basic
property 
\begin{equation}
\underset{g}{\tau}\underset{g^{-1}}{\cdot}\underset{g}{\tau}=\tau \underset{%
g^{-1}}{\lrcorner}\underset{g}{\widetilde{\tau}}=\underset{g}{\tau }\underset%
{g^{-1}}{}\underset{g}{\widetilde{\tau}}=(-1)^{q}.  \label{OHD.7a}
\end{equation}
In order to prove it we should recall that $sgn(\det[g])=(-1)^{q},$ where $q$
is the \emph{number of negative eigenvalues} of $g.$

An expansion formula for smooth pseudoscalar fields on $U$ can be also
obtained from Eq.(\ref{OHD.7a}), i.e., 
\begin{equation}
I=(-1)^{q}(I\underset{g^{-1}}{\cdot}\underset{g}{\tau})\underset{g}{\tau}.
\label{OHD.7b}
\end{equation}

Associated to $\tau,$ the smooth extensor field on $U,$ namely $\star,$
defined by 
\begin{equation}
\star X=\widetilde{X}\lrcorner\tau=\widetilde{X}\tau,  \label{OHD.8}
\end{equation}
will be called the \emph{standard Hodge extensor field }on $U$.

Such $\star$ is \emph{non-singular} and its inverse $\star^{-1}$ is given by 
\begin{equation}
\star^{-1}X=\tau\llcorner\widetilde{X}=\tau\widetilde{X}.  \label{OHD.8a}
\end{equation}

By using a property analogous to that result given by Eq.(\ref{OHD.6b}) we
can easily prove that the standard Hodge extensor field is \emph{ordinarily
constant}, i.e., 
\begin{equation}
a\cdot\partial_{o}\star=0.  \label{OHD.8b}
\end{equation}
Also, 
\begin{equation}
a\cdot\partial_{o}\star^{-1}=0.  \label{OHD.8c}
\end{equation}

Associated to $\underset{g}{\tau},$ the smooth extensor field on $U,$ namely 
$\underset{g}{\star},$ defined by 
\begin{equation}
\underset{g}{\star}X=\widetilde{X}\underset{g^{-1}}{\lrcorner}\underset{g}{%
\tau}=\widetilde{X}\underset{g^{-1}}{}\underset{g}{\tau}=\sqrt{\left| \det[g]%
\right| }\underline{g}^{-1}(\widetilde{X})\lrcorner\tau=\sqrt{\left| \det[g]%
\right| }\underline{g}^{-1}(\widetilde{X})\tau,  \label{OHD.9}
\end{equation}
will be called the \emph{metric Hodge extensor field }on $U.$ Of course, it
is associated to the metric structure $(U,g).$

Such $\underset{g}{\star}$ is also non-singular and its inverse, namely $%
\underset{g}{\star}^{-1},$ is given by 
\begin{align}
\underset{g}{\star}^{-1}X & =(-1)^{q}\underset{g}{\tau}\underset{g^{-1}}{%
\llcorner}\widetilde{X}=(-1)^{q}\underset{g}{\tau}\underset{g^{-1}}{}%
\widetilde{X},  \notag \\
& =(-1)^{q}\sqrt{\left| \det[g]\right| }\tau\llcorner\underline{g}^{-1}(%
\widetilde{X})=(-1)^{q}\sqrt{\left| \det[g]\right| }\tau\underline {g}^{-1}(%
\widetilde{X}).  \label{OHD.9a}
\end{align}

\subsection{Duality Identities}

We present in this subsection two interesting and useful formulas which
relate the \emph{ordinary curl} $\partial_{o}\wedge$ to the \emph{ordinary
contracted divergence} $\partial_{o}\lrcorner.\vspace{0.1in}$

\textbf{i.} For all $X\in\mathcal{M}(U)$ it holds 
\begin{equation}
\tau(\partial_{o}\wedge X)=(-1)^{n+1}\partial_{o}\lrcorner(\tau X).
\label{DI.1}
\end{equation}

To prove Eq.(\ref{DI.1}) we will use the so-called duality identity $%
I(a\wedge Y)=(-1)^{n+1}a\lrcorner(IY),$ where $a\in\mathcal{U}_{o},$ $%
I\in\bigwedge ^{n}\mathcal{U}_{o}$ and $Y\in\bigwedge\mathcal{U}_{o}.$ By
recalling the known identities $\partial_{a}\wedge(a\cdot\partial_{o}X)=%
\partial_{o}\wedge X$ and $\partial_{a}\lrcorner(a\cdot\partial_{o}X)=%
\partial_{o}\lrcorner X,$ and using Eq.(\ref{OHD.6b}), we get 
\begin{align*}
\tau(\partial_{o}\wedge X) & =\tau(\partial_{a}\wedge(a\cdot\partial
_{o}X))=(-1)^{n+1}\partial_{a}\lrcorner(\tau(a\cdot\partial_{o}X)), \\
& =(-1)^{n+1}\partial_{a}\lrcorner(a\cdot\partial_{o}(\tau
X))=(-1)^{n+1}\partial_{o}\lrcorner(\tau X).
\end{align*}

\textbf{ii.} For all $X\in\mathcal{M}(U)$ it holds 
\begin{equation}
\tau\underset{g^{-1}}{}(\partial_{o}\wedge X)=\frac{(-1)^{n+1}}{\det [g]}%
\underline{g}(\partial_{o}\lrcorner(\tau X)).  \label{DI.2}
\end{equation}

We prove (\ref{DI.2}) using the identity\footnote{%
In order to prove it we should use the following identities: $I\underset{%
g^{-1}}{}Y=I\underline {g}^{-1}(Y),$ $X\llcorner\underline{t}^{-1}(Y)=%
\underline{t}^{\dagger }(\underline{t}^{*}(X)\llcorner Y)$ and $\underline{t}%
^{-1}(I)=\left. \det\right. ^{-1}[t]I,$ and so forth.} $I\underset{g^{-1}}{}%
Y=\left. \det\right. ^{-1}[g]\underline{g}(IY)$ and Eq.(\ref{DI.1}).

\subsection{Hodge Duality Identities}

Now, we present two noticeable identities which relate the curl $\partial
_{o}\wedge$ to the contracted divergence $\partial_{o}\lrcorner$ involving
the standard and the metric Hodge extensor fields $\star$ and $\underset{g}{%
\star }$.\vspace{0.1in}

\textbf{i.} For all $X\in\mathcal{M}(U)$ it holds 
\begin{equation}
\star^{-1}(\partial_{o}\wedge(\star X))=-\partial_{o}\lrcorner\widehat{X}.
\label{HDI.1}
\end{equation}

To show Eq.(\ref{HDI.1}) we will use the duality identity given by Eq.(\ref%
{DI.1}). By using Eq.(\ref{OHD.8a}) and Eq.(\ref{OHD.8}), and recalling the
identities $\widetilde{\partial_{o}\wedge Y}=\partial_{o}\wedge\overline{Y}$
and $\overline{XY}=\overline{Y}$ $\overline{X},$ and the obvious property $%
\tau\overline{\tau}=(-1)^{n},$ we get 
\begin{align*}
\star^{-1}(\partial_{o}\wedge(\star X)) & =\tau\widetilde{(\partial
_{o}\wedge(\widetilde{X}\tau))}=\tau(\partial_{o}\wedge(\overline{\tau }%
\widehat{X})), \\
& =(-1)^{n+1}\partial_{o}\lrcorner(\tau\overline{\tau}\widehat{X}%
)=-\partial_{o}\lrcorner\widehat{X}.
\end{align*}

\textbf{ii.} For all $X\in\mathcal{M}(U)$ it holds 
\begin{equation}
\underset{g}{\star}^{-1}(\partial_{o}\wedge(\underset{g}{\star}X))=-\frac {1%
}{\sqrt{\left| \det[g]\right| }}\underline{g}(\partial_{o}\lrcorner (\sqrt{%
\left| \det[g]\right| }\underline{g}^{-1}(\widehat{X}))).  \label{HDI.2}
\end{equation}

This follows from appropriate use of the duality identity given by Eq.(\ref%
{DI.2}). Indeed, a straightforward calculation using Eq.(\ref{OHD.9a}) and
Eq.(\ref{OHD.9}) allows us to write 
\begin{align*}
\underset{g}{\star}^{-1}(\partial_{o}\wedge(\underset{g}{\star}X)) &
=(-1)^{q}\sqrt{\left\vert \det[g]\right\vert }\tau\underset{g^{-1}}{}%
\widetilde{(\partial_{o}\wedge(\widetilde{X}\underset{g^{-1}}{}\underset{g}{%
\tau}))} \\
& =(-1)^{q}\sqrt{\left\vert \det[g]\right\vert }\tau\underset{g^{-1}}{}%
\partial_{o}\wedge(\sqrt{\left\vert \det[g]\right\vert }\overline{\tau }%
\underline{g}^{-1}(\widehat{X})) \\
& =(-1)^{q}\sqrt{\left\vert \det[g]\right\vert }\frac{(-1)^{n+1}}{\det [g]}%
\underline{g}(\partial_{o}\lrcorner(\sqrt{\left\vert \det[g]\right\vert }\tau%
\overline{\tau}\underline{g}^{-1}(\widehat{X}))), \\
& =-\frac{1}{\sqrt{\left\vert \det[g]\right\vert }}\underline{g}(\partial
_{o}\lrcorner(\sqrt{\left\vert \det[g]\right\vert }\underline{g}^{-1}(%
\widehat{X}))).
\end{align*}

\subsection{Ordinary Hodge Coderivative Operators}

We introduce the so-called \emph{standard Hodge derivative operator} $\delta:%
\mathcal{M}(U)\rightarrow\mathcal{M}(U)$ such that 
\begin{equation}
\delta X=\star^{-1}(\partial_{o}\wedge(\star\widehat{X})).  \label{OHO.1}
\end{equation}

Its basic property is 
\begin{equation}
\delta X=-\partial_{o}\lrcorner X.  \label{OHO.1a}
\end{equation}
It immediately follows from the Hodge duality identity given by Eq.(\ref%
{HDI.1}).

We introduce the so-called \emph{metric Hodge coderivative operator} $%
\underset{g}{\delta}:\mathcal{M}(U)\rightarrow\mathcal{M}(U)$ such that 
\begin{equation}
\underset{g}{\delta}X=\text{ }\underset{g}{\star}^{-1}(\partial_{o}\wedge(%
\underset{g}{\star}\widehat{X})).  \label{OHO.2}
\end{equation}

It satisfies the basic property 
\begin{equation}
\underset{g}{\delta}X=-\frac{1}{\sqrt{\left\vert \det[g]\right\vert }}%
\underline{g}(\partial_{o}\lrcorner(\sqrt{\left\vert \det[g]\right\vert }%
\underline{g}^{-1}(X))),  \label{OHO.2a}
\end{equation}
which is an immediate consequence of the Hodge duality identity given by Eq.(%
\ref{HDI.2}).

\subsection{ Levi-Civita Geometric Structure}

We recall from Section 2 that the \emph{Levi-Civita connection field} is the
smooth vector elementary $2$-extensor field on $U,$ namely $\lambda ,$
defined by 
\begin{equation}
\lambda (a,b)=\frac{1}{2}g^{-1}\circ (a\cdot \partial _{o}g)(b)+\omega
_{0}(a)\underset{g}{\times }b,  \label{LGS.1}
\end{equation}%
where $\omega _{0}$ is the smooth $(1,2)$-extensor field on $U$ given by 
\begin{equation}
\omega _{0}(a)=-\frac{1}{4}\underline{g}^{-1}(\partial _{b}\wedge \partial
_{c})a\cdot ((b\cdot \partial _{o}g)(c)-(c\cdot \partial g_{o})(b)).
\label{LGS.1a}
\end{equation}

Such $\omega_{0}$ satisfies 
\begin{equation}
\omega_{0}(a)\underset{g}{\times}b\underset{g}{\cdot}c=\frac{1}{2}%
a\cdot((b\cdot\partial_{o}g)(c)-(c\cdot\partial_{o}g)(b)).  \label{LGS.1b}
\end{equation}

The \emph{open set} $U$ endowed with $\lambda $ and $g,$ namely $(U,\lambda
,g),$ is a \emph{geometric structure} on $U$, a statement that means that 
\emph{Levi-Civita parallelism structure} $(U,\lambda )$ is compatible with
the metric structure $(U,g).$ Or equivalently, the pair of $a$-\emph{DCDO's }%
associated to $(U,\lambda ),$ namely $(D_{a}^{+},D_{a}^{-}),$ is $g$%
-compatible.

The \emph{Levi-Civita} $a$-\emph{DCDO's }$D_{a}^{+}$ and $D_{a}^{-}$ are
defined by 
\begin{align}
D_{a}^{+}X & =a\cdot\partial_{o}X+\Lambda_{a}(X),  \label{LGS.2a} \\
D_{a}^{-}X & =a\cdot\partial_{o}X-\Lambda_{a}^{\dagger}(X).  \label{LGS.2b}
\end{align}
Note that $\Lambda_{a}$ is the so-called \emph{generalized} of $\lambda_{a}.$
The latter is the so-called $a$-\emph{directional connection field}
associated to $\lambda,$ given by $\lambda_{a}(b)=\lambda(a,b).$

We present now two pairs of noticeable properties of $\lambda.\vspace{0.1in} 
$

\textbf{i. }The \emph{scalar divergence} of $\lambda_{a}(b)$ with respect to
the first variable, namely $\partial_{a}\cdot\lambda_{a}(b),$ and the \emph{%
curl} of $\lambda_{a}^{\dagger}(b)$ with respect to the first variable,
namely $\partial_{a}\wedge\lambda_{a}^{\dagger}(b),$ are given by 
\begin{align}
\partial_{a}\cdot\lambda_{a}(b) & =\frac{1}{\sqrt{\left| \det[g]\right| }}%
b\cdot\partial_{o}\sqrt{\left| \det[g]\right| },  \label{LGS.3a} \\
\partial_{a}\wedge\lambda_{a}^{\dagger}(b) & =0.  \label{LGS.3b}
\end{align}

In order to prove the first result we will use the formula $\tau^{\ast
}(\partial_{n})\cdot(a\cdot\partial_{o}\tau)(n)=\left. \det\right.
^{-1}[\tau]a\cdot\partial_{o}\det[\tau]$, valid for all non-singular smooth $%
(1,1)$-extensor field $\tau.$ By using the symmetry property $%
(g^{-1})^{\dagger}=g^{-1}$ and Eq.(\ref{LGS.1b}), we can write 
\begin{align*}
\partial_{a}\cdot\lambda_{a}(b) & =\frac{1}{2}\partial_{a}\cdot(g^{-1}%
\circ(a\cdot\partial_{o}g)(b)+\partial_{a}\cdot(\omega_{0}(a)\underset{g}{%
\times}b) \\
& =\frac{1}{2}g^{-1}(\partial_{a})\cdot(a\cdot\partial_{o}g)(b)+\omega
_{0}(g^{-1}(\partial_{a}))\underset{g}{\times}b\underset{g}{\cdot}a \\
& =\frac{1}{2}g^{-1}(\partial_{a})\cdot(b\cdot\partial_{o}g)(a)=\frac{1}{2}%
\frac{1}{\det[g]}b\cdot\partial_{o}\det[g].
\end{align*}
Then, recalling the identity $(f)^{-1}b\cdot\partial_{o}f=2(\left\vert
f\right\vert )^{-1/2}b\cdot\partial_{o}(\left\vert f\right\vert )^{1/2}$
valid for all \emph{non-zero} $f\in\mathcal{S}(U)$, the expected result
immediately follows.

To prove the second result we only need to take into account the \emph{%
symmetry property} $\lambda_{a}(b)=\lambda_{b}(a).$ We have 
\begin{equation*}
\partial_{a}\wedge\lambda_{a}^{\dagger}(b)=\partial_{a}\wedge\partial
_{n}(n\cdot\lambda_{a}^{\dagger}(b))=\partial_{a}\wedge\partial_{n}(%
\lambda_{a}(n)\cdot b)=0.
\end{equation*}

\textbf{ii. }The \emph{left contracted divergence} of $\Lambda_{a}(X)$ with
respect to $a,$ namely $\partial_{a}\lrcorner\Lambda_{a}(X),$ and the \emph{%
curl} of $\Lambda_{a}^{\dagger}(X)$ with respect to $a,$ namely $%
\partial_{a}\wedge\Lambda_{a}^{\dagger}(X),$ are given by 
\begin{align}
\partial_{a}\lrcorner\Lambda_{a}(X) & =\frac{1}{\sqrt{\left| \det [g]\right| 
}}(\partial_{o}\sqrt{\left| \det[g]\right| })\lrcorner X,  \label{LGS.4a} \\
\partial_{a}\wedge\Lambda_{a}^{\dagger}(X) & =0.  \label{LGS.4b}
\end{align}

In order to prove Eq.(\ref{LGS.4a}) we will use the multivector identities $%
v\lrcorner(X\wedge Y)=(v\lrcorner X)\wedge Y+\widehat{X}\wedge(v\lrcorner Y)$
and $X\lrcorner(Y\lrcorner Z)=(X\wedge Y)\lrcorner Z,$ where $v\in \mathcal{U%
}_{o}$ and $X,Y,Z\in\bigwedge\mathcal{U}_{o}.$ A straightforward calculation
using Eq.(\ref{LGS.3a}) allows us to get 
\begin{align*}
\partial_{a}\lrcorner\Lambda_{a}(X) &
=\partial_{a}\cdot\lambda_{a}(\partial_{n})(n\lrcorner
X)-\lambda_{a}(\partial_{n})\wedge(\partial _{a}\lrcorner(n\lrcorner X)) \\
& =\frac{1}{\sqrt{\left\vert \det[g]\right\vert }}(\partial_{n}\cdot
\partial_{o}\sqrt{\left\vert \det[g]\right\vert })(n\lrcorner X)-\lambda
_{a}(\partial_{n})\wedge((\partial_{a}\wedge n)\lrcorner X), \\
& =\frac{1}{\sqrt{\left\vert \det[g]\right\vert }}\partial_{n}(n\cdot
\partial_{o}\sqrt{\left\vert \det[g]\right\vert })\lrcorner X-\lambda
_{\partial_{a}}(\partial_{n})\wedge((a\wedge n)\lrcorner X).
\end{align*}
Then, recalling the identity $\partial_{n}(n\cdot\partial_{o}Y)=%
\partial_{o}Y $ and the symmetry property $\lambda_{a}(b)=\lambda_{b}(a),$
we get the required result.

The proof of the second property follows immediately by using Eq.(\ref%
{LGS.3b}), and recalling that the adjoint of generalized is equal to the
generalized of adjoint, i.e., $\Lambda_{a}^{\dagger}(X)=\lambda_{a}^{%
\dagger}(\partial _{n})\wedge(n\lrcorner X)$.

\subsection{Levi-Civita Derivatives}

We introduce now the \emph{canonical covariant divergence operator}, $%
D^{+}\lrcorner:\mathcal{M}(U)\rightarrow\mathcal{M}(U)$ such that 
\begin{equation}
D^{+}\lrcorner X=\partial_{a}\lrcorner(D_{a}^{+}X),  \label{LCD.1}
\end{equation}
i.e., $D^{+}\lrcorner
X=e^{\mu}\lrcorner(D_{e_{\mu}}^{+}X)=e_{\mu}\lrcorner(D_{e^{\mu}}^{+}X),$
where $(\{e_{\mu}\},\{e^{\mu}\})$ is any pair of \emph{canonical reciprocal
frame fields} on $U.$

Its basic property is 
\begin{equation}
D^{+}\lrcorner X=\frac{1}{\sqrt{\left\vert \det[g]\right\vert }}\partial
_{o}\lrcorner(\sqrt{\left\vert \det[g]\right\vert }X).  \label{LCD.1a}
\end{equation}
Indeed, using the identity $\partial_{a}\lrcorner(a\cdot\partial
_{o}X)=\partial_{o}\lrcorner X$ and Eq.(\ref{LGS.4a}) into the definition
given by Eq.(\ref{LGS.2a}), we get 
\begin{align*}
D^{+}\lrcorner X & =\partial_{a}\lrcorner(a\cdot\partial_{o}X)+\partial
_{a}\lrcorner\Lambda_{a}(X), \\
& =\partial\lrcorner X+\frac{1}{\sqrt{\left\vert \det[g]\right\vert }}%
(\partial_{o}\sqrt{\left\vert \det[g]\right\vert })\lrcorner X.
\end{align*}
So, by recalling the identity $\partial_{o}\lrcorner(fY)=(\partial
_{o}f)\lrcorner X+f(\partial_{o}\lrcorner X),$ for all $f\in\mathcal{S}(U)$
and $Y\in\mathcal{M}(U),$ we can get the proof for this remarkable property.

The so-called \emph{metric covariant divergence}, \emph{covariant curl }and 
\emph{metric covariant gradient operators}, namely $D^{-}\underset{g^{-1}}{%
\lrcorner},$ $D^{-}\wedge$ and $D^{-}\underset{g^{-1}}{},$ map also smooth
multivector fields to smooth multivector fields and are defined by 
\begin{align}
D^{-}\underset{g^{-1}}{\lrcorner}X & =\partial_{a}\underset{g^{-1}}{\lrcorner%
}(D_{a}^{-}X)=g^{-1}(\partial_{a})\lrcorner(D_{a}^{-}X),  \label{LCD.2a} \\
D^{-}\wedge X & =\partial_{a}\wedge(D_{a}^{-}X),  \label{LCD.2b} \\
D^{-}\underset{g^{-1}}{}X & =\partial_{a}\underset{g^{-1}}{}(D_{a}^{-}X).
\label{LCD.2c}
\end{align}
These operators are related by 
\begin{equation}
D^{-}\underset{g^{-1}}{}X=D^{-}\underset{g^{-1}}{\lrcorner}X+D^{-}\wedge X.
\label{LCD.3}
\end{equation}

The basic properties of the metric covariant divergence are 
\begin{align}
D^{-}\underset{g^{-1}}{\lrcorner}X & =\underline{g}(D^{+}\lrcorner 
\underline{g}^{-1}(X)),  \label{LCD.4a} \\
D^{-}\underset{g^{-1}}{\lrcorner}X & =\frac{1}{\sqrt{\left| \det[g]\right| }}%
\underline{g}(\partial_{o}\lrcorner(\sqrt{\left| \det[g]\right| }\underline{g%
}^{-1}(X))).  \label{LCD.4b}
\end{align}

Eq.(\ref{LCD.4a}) follows from the fundamental property $D_{a}^{-}\underline{%
g}(X)=\underline{g}(D_{a}^{+}X)$ which holds for any $g$-compatible pair of $%
a$-\emph{DCDO's, }by using the identity $X\lrcorner\underline {t}(Y)=%
\underline{t}(\underline{t}^{\dagger}(X)\lrcorner Y).$\emph{\ }Eq.(\ref%
{LCD.4b}) is deduced by using Eq.(\ref{LCD.4a}) and Eq.(\ref{LCD.1a}).

A remarkable property which follows from Eq.(\ref{LCD.4b}) is 
\begin{equation}
D^{-}\underset{g^{-1}}{\lrcorner}(D^{-}\underset{g^{-1}}{\lrcorner}X)=0.
\label{LCD.4b1}
\end{equation}
In order to prove it we should use the known identity $\partial_{o}%
\lrcorner(\partial_{o}\lrcorner Y)=0.$

By comparing Eq.(\ref{OHO.2a}) and Eq.(\ref{LCD.4b}) we get 
\begin{equation}
D^{-}\underset{g^{-1}}{\lrcorner}X=-\underset{g}{\delta}X.  \label{LCD.4b2}
\end{equation}

The basic property for the covariant curl is 
\begin{equation}
D^{-}\wedge X=\partial_{o}\wedge X.  \label{LCD.5}
\end{equation}
It immediately follows by using the identity $\partial_{a}\wedge
(a\cdot\partial_{o}X)=\partial_{o}\wedge X$ and Eq.(\ref{LGS.4b}) into the
definition given by Eq.(\ref{LGS.2b}).

The four \emph{derivative-like operators} defined by Eq.(\ref{LCD.1}), Eq.(%
\ref{LCD.2a}), Eq.(\ref{LCD.2b}) and Eq.(\ref{LCD.2c}) will be called the 
\emph{Levi-Civita derivatives. }They are involved in three useful identities
which are used in the \emph{Lagrangian theory of multivector fields }(\cite%
{3,rodoliv2006}). These are, 
\begin{align}
& (\partial_{o}\wedge X)\underset{g^{-1}}{\cdot}Y+X\underset{g^{-1}}{\cdot }%
(D^{-}\underset{g^{-1}}{\lrcorner}Y)  \notag \\
& =\frac{1}{\sqrt{\left\vert \det[g]\right\vert }}\partial_{o}\cdot (\sqrt{%
\left\vert \det[g]\right\vert }\partial_{n}(n\wedge X)\underset{g^{-1}}{\cdot%
}Y),  \label{LCD.6a} \\
& (D^{-}\underset{g^{-1}}{\lrcorner}X)\underset{g^{-1}}{\cdot}Y+X\underset{%
g^{-1}}{\cdot}(\partial_{o}\wedge Y)  \notag \\
& =\frac{1}{\sqrt{\left\vert \det[g]\right\vert }}\partial_{o}\cdot (\sqrt{%
\left\vert \det[g]\right\vert }\partial_{n}(n\underset{g^{-1}}{\lrcorner}X)%
\underset{g^{-1}}{\cdot}Y),  \label{LCD.6b} \\
& (D^{-}\underset{g^{-1}}{}X)\underset{g^{-1}}{\cdot}Y+X\underset{g^{-1}}{%
\cdot}(D^{-}\underset{g^{-1}}{}Y)  \notag \\
& =\frac{1}{\sqrt{\left\vert \det[g]\right\vert }}\partial_{o}\cdot (\sqrt{%
\left\vert \det[g]\right\vert }\partial_{n}(n\underset{g^{-1}}{}X)\underset{%
g^{-1}}{\cdot}Y).  \label{LCD.6c}
\end{align}

\subsection{Gauge (Deformed) Derivatives}

Let $h$ be a gauge metric field for $g$. This statement, as we already know
means that there is a smooth $(1,1)$-extensor field $h$ such that $%
g=h^{\dagger }\circ \eta \circ h,$ where $\eta $ is an orthogonal metric
field with the same signature as $g.$ As we know, associated to a
Levi-Civita $a$-\emph{DCDO's,} namely $(D_{a}^{+},D_{a}^{-}),$ there must be
an unique $\eta $-compatible pair of $a$-\emph{DCDO's, }namely\emph{\ }$(%
\mathbf{D}_{a}^{+},\mathbf{D}_{a}^{-}),$ given by the following formulas 
\begin{align}
\text{$\mathbf{D}$}_{ha}^{+}X& =\underline{h}(D_{a}^{+}\underline{h}%
^{-1}(X)),  \label{GD.1a} \\
\text{$\mathbf{D}$}_{h^{\ast }a}^{-}X& =\underline{h}^{\ast }(D_{a}^{+}%
\underline{h}^{\dagger }(X)).  \label{GD.1b}
\end{align}%
These equations say that $(\mathbf{D}_{ha}^{+},\mathbf{D}_{h^{\ast }a}^{-})$
is the $h$-\emph{deformation} of $(D_{a}^{+},D_{a}^{-}).$ So, $\mathbf{D}%
_{ha}^{+}$ and $\mathbf{D}$\ $_{h^{\ast }a}^{-}$ will be called the \emph{%
gauge (deformed) covariant derivatives} associated to $D_{a}^{+}$ and $%
D_{a}^{-}.$

We present here two noticeable properties of $\mathbf{D}_{ha}^{+}.\vspace{%
0.1in}$

\textbf{i.} For all $a,b,c\in \mathcal{V}(U),$ it holds 
\begin{equation}
(\text{$\mathbf{D}$}_{ha}^{+}b)\underset{\eta }{\cdot }%
c=[h(a),h^{-1}(b),h^{-1}(c)].  \label{GD.2}
\end{equation}

To prove Eq.(\ref{GD.2}) we use Eq.(\ref{GD.1a}), the identity $%
(D_{a}^{+}b)\cdot c=\QDATOPD\{ \} {c}{a,b},$ and the definition of the
Christoffel operator of second kind, i.e., $\QDATOPD\{ \}
{c}{a,b}=[a,b,g^{-1}(c)]$. Indeed, we can write 
\begin{align*}
(\text{$\mathbf{D}$}_{ha}^{+}b)\underset{\eta }{\cdot }c&
=D_{a}^{+}(h^{-1}(b))\cdot h^{\dagger }\circ \eta (c) \\
& =\QDATOPD\{ \} {h^{\dagger }\circ \eta
(c)}{h(a),h^{-1}(b)}=[h(a),h^{-1}(b),g^{-1}\circ h^{\dagger }\circ \eta (c)],
\end{align*}%
and recalling that $g^{-1}=h^{-1}\circ \eta \circ h^{\ast },$ the required
result immediately follows.

\textbf{ii.} There exists a smooth $(1,2)$-extensor field on $U_{o}$, namely 
$\Omega _{0},$ such that 
\begin{equation}
\text{$\mathbf{D}$}_{a}^{+}X=a\cdot \partial _{o}X+\Omega _{0}(a)\underset{%
\eta }{\times }X.  \label{GD.3}
\end{equation}%
Such $\Omega _{0}$ is given by 
\begin{equation}
\Omega _{0}(a)=-\frac{1}{2}\underline{\eta }(\partial _{b}\wedge \partial
_{c})[a,h^{-1}(b),h^{-1}(c)].  \label{GD.4}
\end{equation}

The proof of this important result is as follows. First we prove a
particular case of the above property, i.e., 
\begin{equation}
\text{$\mathbf{D}$}_{a}^{+}b=a\cdot \partial _{o}b+\Omega _{0}(a)\underset{%
\eta }{\times }b.  \label{GD.3a}
\end{equation}

To do so, we will use the following properties of the Christoffel operator
of first kind $[a,b,c]-[a,c,b]=2[a,b,c]-a\cdot \partial _{o}(b\underset{g}{%
\cdot }c),$ and $[a,b+b^{\prime },c]=[a,b,c]+[a,b^{\prime },c]$ and $%
[a,fb,c]=f[a,b,c]+(a\cdot \partial _{o}f)b\underset{g}{\cdot }c.$ A
straightforward calculation permit us to write 
\begin{align*}
\Omega _{0}(a)\underset{\eta }{\times }b\underset{\eta }{\cdot }c& =(b\wedge
c)\underset{\eta }{\cdot }\Omega _{0}(a)=\frac{1}{2}(b\wedge c)\cdot
(\partial _{p}\wedge \partial _{q})[a,h^{-1}(p),h^{-1}(q)] \\
& =\frac{1}{2}\det \left[ 
\begin{array}{cc}
b\cdot \partial _{p} & b\cdot \partial _{q} \\ 
c\cdot \partial _{p} & c\cdot \partial _{q}%
\end{array}%
\right] [a,h^{-1}(p),h^{-1}(q)] \\
& =\frac{1}{2}b\cdot \partial _{p}c\cdot \partial \lbrack
a,h^{-1}(p),h^{-1}(q)]-b\cdot \partial _{q}c\cdot \partial
_{p}[a,h^{-1}(p),h^{-1}(q)] \\
& =\frac{1}{2}b\cdot \partial _{p}c\cdot \partial
_{q}([a,h^{-1}(p),h^{-1}(q)]-[a,h^{-1}(q),h^{-1}(p)]) \\
& =\frac{1}{2}b\cdot \partial _{p}c\cdot \partial
_{q}(2[a,h^{-1}(p),h^{-1}(q)]-a\cdot \partial _{o}(h^{-1}(p)\underset{g}{%
\cdot }h^{-1}(q))) \\
& =b\cdot \partial _{p}c\cdot \partial _{q}(p\cdot b_{\mu }[a,h^{-1}(b_{\mu
}),h^{-1}(c)]+(a\cdot \partial _{o}p)\underset{\eta }{\cdot }c) \\
& =b\cdot b_{\mu }[a,h^{-1}(b_{\mu }),h^{-1}(c)]+(a\cdot \partial _{o}b)%
\underset{\eta }{\cdot }c-(a\cdot \partial _{o}b)\underset{\eta }{\cdot }c \\
& =[a,h^{-1}(b),h^{-1}(c)]-(a\cdot \partial _{o}b)\underset{\eta }{\cdot }c,
\end{align*}%
hence, by using Eq.(\ref{GD.2}), the particular case given by Eq.(\ref{GD.3a}%
) immediately follows. Now, we can prove the general case of the above
property .

As we can see, the $a$-directional connection field on $U$ for $(\mathbf{D}%
_{a}^{+},\mathbf{D}_{a}^{-})$ is given by $b\mapsto \Omega _{0}(a)\underset{%
\eta }{\times }b.$ Then, its generalized (extensor field) must be given by $%
X\mapsto (\Omega _{0}(a)\underset{\eta }{\times }\partial _{b})\wedge
(b\lrcorner X).$ But, by recalling the noticeable identity $(B\underset{\eta 
}{\times }\partial _{b})\wedge (b\lrcorner X)=B\underset{\eta }{\times }X,$
where $B\in \bigwedge^{2}\mathcal{U}_{o}$ and $X\in \bigwedge \mathcal{U}%
_{o},$ we find that it can be written as $X\mapsto \Omega _{0}(a)\underset{%
\eta }{\times }X$.

We introduce now the \emph{gauge covariant divergence}, \emph{gauge
covariant curl} and \emph{gauge covariant gradient operators}, namely $%
\mathbf{D}^{-}\underset{\eta }{\lrcorner },$ $\mathbf{D}^{-}\wedge $ and $%
\mathbf{D}^{-}\underset{\eta }{}.$ They all map smooth multivector fields to
smooth multivector fields and are defined by 
\begin{align}
\text{$\mathbf{D}$}^{-}\underset{\eta }{\lrcorner }X& =h^{\ast }(\partial
_{a})\underset{\eta }{\lrcorner }(\text{$\mathbf{D}$}_{h^{\ast }a}^{-}X),
\label{GD.5a} \\
\text{$\mathbf{D}$}^{-}\wedge X& =h^{\ast }(\partial _{a})\wedge (\text{$%
\mathbf{D}$}_{h^{\ast }a}^{-}X),  \label{GD.5b} \\
\text{$\mathbf{D}$}^{-}\underset{\eta }{}X& =h^{\ast }(\partial _{a})%
\underset{\eta }{}(\text{$\mathbf{D}$}_{h^{\ast }a}^{-}X).  \label{GD.5c}
\end{align}%
It is obvious that the relationship among them is given by 
\begin{equation}
\text{$\mathbf{D}$}^{-}\underset{\eta }{}X=\text{$\mathbf{D}$}^{-}\underset{%
\eta }{\lrcorner }X+\text{$\mathbf{D}$}^{-}\wedge X.  \label{GD.6}
\end{equation}

Their basic properties are given by the following \emph{golden formulas} 
\begin{align}
\underline{h}^{\ast }(D^{-}\underset{g^{-1}}{\lrcorner }X)& =\text{$\mathbf{D%
}$}^{-}\underset{\eta }{\lrcorner }\underline{h}^{\ast }(X),  \label{GD.7a}
\\
\underline{h}^{\ast }(D^{-}\wedge X)& =\text{$\mathbf{D}$}^{-}\wedge 
\underline{h}^{\ast }(X),  \label{GD.7b} \\
\underline{h}^{\ast }(D^{-}\underset{g^{-1}}{}X)& =\text{$\mathbf{D}$}^{-}%
\underset{\eta }{}\underline{h}^{\ast }(X).  \label{GD.7c}
\end{align}%
These formulas can be proved by using the golden formula deduced in \cite{2}%
, $\underline{h}^{\ast }(X\underset{g^{-1}}{\ast }Y)=\underline{h}^{\ast }(X)%
\underset{\eta }{\ast }\underline{h}^{\ast }(Y),$ where $X,Y\in \bigwedge 
\mathcal{U}_{o},$ and $\underset{g^{-1}}{\ast }$ means either exterior
product or any $g^{-1}$-product of smooth multivector fields, and
analogously for $\underset{\eta }{\ast }$. It is also necessary to take into
account the master formulas $g=h^{\dagger }\circ \eta \circ h$ and $%
g^{-1}=h^{-1}\circ \eta \circ h^{\ast },$ and the relationship between $%
\mathbf{D}_{ha}^{+}$ and $\mathbf{D}_{h^{\ast }a}^{-},$ i.e., $\underline{%
\eta }(\mathbf{D}_{ha}^{+}\underline{\eta }(X))=$ $\mathbf{D}_{h^{\ast
}a}^{-}X.$

From Eq.(\ref{LCD.4b}) and Eq.(\ref{LCD.5}) by using Eq.(\ref{GD.7a}) and
Eq.(\ref{GD.7b}) we find the interesting identities 
\begin{align}
\text{$\mathbf{D}$}^{-}\underset{\eta }{\lrcorner }X& =\frac{1}{\det [h]}%
\underline{\eta \circ h}(\partial _{o}\lrcorner (\det [h]\underline{%
h^{-1}\circ \eta }(X))),  \label{GD.8a} \\
\text{$\mathbf{D}$}^{-}\wedge X& =\underline{h}^{\ast }(\partial _{o}\wedge 
\underline{h}^{\dagger }(X)).  \label{GD.8b}
\end{align}

\subsection{Covariant Hodge Coderivative}

Let $(U,\gamma,g)$ be a \emph{geometric structure} on $U,$ and let us denote
by $(\mathcal{D}_{a}^{+},\mathcal{D}_{a}^{-})$ the $g$-\emph{compatible}
pair of $a$-\emph{DCDO's }associated to $(U,\gamma,g).$ We will present two
noticeable properties which are satisfied by the second $a$-\emph{DCDO.}

\textbf{i.} $\underset{g}{\tau}$ is a \emph{covariantly constant} smooth
pseudoscalar field on $U,$ i.e., 
\begin{equation}
\mathcal{D}_{a}^{-}\underset{g}{\tau}=0.  \label{CHC.1}
\end{equation}

To show Eq.(\ref{CHC.1}) we will use the basic property and the expansion
formula given by Eq.(\ref{OHD.7a}) and Eq.(\ref{OHD.7b}). From the
Ricci-like theorem for $\mathcal{D}_{a}^{-},$ we have that 
\begin{equation*}
(\mathcal{D}_{a}^{-}\underset{g}{\tau})\underset{g^{-1}}{\cdot}\underset{g}{%
\tau}+\underset{g}{\tau}\underset{g^{-1}}{\cdot}(\mathcal{D}_{a}^{-}\underset%
{g}{\tau})=0,
\end{equation*}
i.e., $(\mathcal{D}_{a}^{-}\underset{g}{\tau})\underset{g^{-1}}{\cdot }%
\underset{g}{\tau}=0.$ Then, $\mathcal{D}_{a}^{-}\underset{g}{\tau}%
=(-1)^{q}((\mathcal{D}_{a}^{-}\underset{g}{\tau})\underset{g^{-1}}{\cdot }%
\underset{g}{\tau})\underset{g}{\tau}=0$.

\textbf{ii.} For all $X\in\mathcal{M}(U)$ it holds 
\begin{equation}
\mathcal{D}_{a}^{-}(\underset{g}{\tau}\underset{g^{-1}}{*}X)=\underset{g}{%
\tau}\underset{g^{-1}}{*}\mathcal{D}_{a}^{-}(X),  \label{CHC.2}
\end{equation}
where $\underset{g^{-1}}{*}$ means either exterior product or any $g^{-1}$%
-product of smooth multivector fields.

This result follows at once using Eq.(\ref{CHC.1}) and the general Leibniz
rule for $\mathcal{D}_{a}^{-}$.

We introduce now the \emph{metric covariant divergence}, \emph{covariant curl%
} and \emph{metric covariant gradient operators}, namely $\mathcal{D}^{-}%
\underset{g^{-1}}{\lrcorner },$ $\mathcal{D}^{-}\wedge $ and $\mathcal{D}^{-}%
\underset{g^{-1}}{},$ which map smooth multivector fields to smooth
multivector fields. They are defined by 
\begin{align}
\mathcal{D}^{-}\underset{g^{-1}}{\lrcorner }X& =\partial _{a}\underset{g^{-1}%
}{\lrcorner }(\mathcal{D}_{a}^{-}X),  \label{CHC.3a} \\
\mathcal{D}^{-}\wedge X& =\partial _{a}\wedge (\mathcal{D}_{a}^{-}X),
\label{CHC.3b} \\
\mathcal{D}^{-}\underset{g^{-1}}{}X& =\partial _{a}\underset{g^{-1}}{}(%
\mathcal{D}_{a}^{-}X).  \label{CHC.3c}
\end{align}%
The relationship among them is given by

\begin{equation}
\mathcal{D}^{-}\underset{g^{-1}}{}X=\mathcal{D}^{-}\underset{g^{-1}}{%
\lrcorner}X+\mathcal{D}^{-}\wedge X.  \label{CHC.3d}
\end{equation}

Now, we will present two noticeable duality identities between $\mathcal{D}%
^{-}\wedge$ and $\mathcal{D}^{-}\underset{g^{-1}}{\lrcorner}.$ One of them
involves the metric Hodge extensor field $\underset{g}{\star}.$

\textbf{iii. }For all $X\in \mathcal{M}(U)$ it holds 
\begin{equation}
\underset{g}{\mathbf{t}}\underset{g^{-1}}{}(\mathcal{D}^{-}\wedge
X)=(-1)^{n+1}\mathcal{D}^{-}\underset{g^{-1}}{\lrcorner }(\underset{g}{%
\mathbf{t}}\underset{g^{-1}}{}X).  \label{CHC.4}
\end{equation}

To show Eq.(\ref{CHC.4}) we will use the $g^{-1}$-duality identity $I%
\underset{g^{-1}}{}(a\wedge Y)=(-1)^{n+1}a\underset{g^{-1}}{\lrcorner }(I%
\underset{g^{-1}}{}Y),$ where $a\in \mathcal{U}_{o},$ $I\in \bigwedge^{n}%
\mathcal{U}_{o}$ and $Y\in \bigwedge \mathcal{U}_{o}.$ A straightforward
calculation by taking into account Eq.(\ref{CHC.2}) gives 
\begin{align*}
\underset{g}{\tau }\underset{g^{-1}}{}(\mathcal{D}^{-}\wedge X)&
=(-1)^{n+1}\partial _{a}\underset{g^{-1}}{\lrcorner }(\underset{g}{\tau }%
\underset{g^{-1}}{}(\mathcal{D}_{a}^{-}X)) \\
& =(-1)^{n+1}\partial _{a}\underset{g^{-1}}{\lrcorner }(\mathcal{D}_{a}^{-}(%
\underset{g}{\tau }\underset{g^{-1}}{}X)), \\
& =(-1)^{n+1}\mathcal{D}^{-}\underset{g^{-1}}{\lrcorner }(\underset{g}{\tau }%
\underset{g^{-1}}{}X).
\end{align*}

\textbf{iv. }For all $X\in\mathcal{M}(U)$ it holds 
\begin{equation}
\underset{g}{\star}^{-1}(\mathcal{D}^{-}\wedge(\underset{g}{\star }X))=-%
\mathcal{D}^{-}\underset{g^{-1}}{\lrcorner}\widehat{X}.  \label{CHC.5}
\end{equation}

This follows by appropriate use of the duality identity given by Eq.(\ref%
{CHC.4}). Indeed, using Eq.(\ref{OHD.9a}) Eq.(\ref{OHD.9}), and recalling
the identities $\widetilde{\mathcal{D}^{-}\wedge X}=\mathcal{D}^{-}\wedge%
\overline{X}$ and $\overline{XY}=\overline{Y}$ $\overline{X},$ and the
obvious property $\underset{g}{\tau}\underset{g^{-1}}{}\overline {\underset{g%
}{\tau}}=(-1)^{n+q}$, we get 
\begin{align*}
\underset{g}{\star}^{-1}(\mathcal{D}^{-}\wedge(\underset{g}{\star}X)) &
=(-1)^{q}\underset{g}{\tau}\underset{g^{-1}}{}\widetilde{(\mathcal{D}%
^{-}\wedge(\widetilde{X}\underset{g^{-1}}{}\underset{g}{\tau}))} \\
& =(-1)^{q}\underset{g}{\tau}\underset{g^{-1}}{}(\mathcal{D}^{-}\wedge(%
\overline{\underset{g}{\tau}}\underset{g^{-1}}{}\widehat{X}) \\
& =(-1)^{q}(-1)^{n+1}\mathcal{D}^{-}\underset{g^{-1}}{\lrcorner}(\underset{g}%
{\tau}\underset{g^{-1}}{}\overline{\underset{g}{\tau}}\underset{g^{-1}}{}%
\widehat{X}), \\
& =-\mathcal{D}^{-}\underset{g^{-1}}{\lrcorner}\widehat{X}.
\end{align*}

We introduce now the \emph{covariant Hodge coderivative operator}, denoted $%
\underset{g}{\Delta}$, and defined by 
\begin{equation*}
\underset{g}{\Delta}:\mathcal{M}(U)\rightarrow\mathcal{M}(U),
\end{equation*}%
\begin{equation}
\underset{g}{\Delta}X=\underset{g}{\text{ }\star}^{-1}(\mathcal{D}^{-}\wedge(%
\underset{g}{\star}\widehat{X})).  \label{CHC.6}
\end{equation}
Its basic property is%
\begin{equation*}
\mathcal{D}^{-}\underset{g^{-1}}{\lrcorner}X=-\underset{g}{\Delta}X,
\end{equation*}
which follows trivially from Eq.(\ref{CHC.5}).

\section{Conclusions}

In this paper we studied how the geometrical calculus of multivector and
extensor fields can be successfully applied to the study of the differential
geometry of an arbitrary smooth manifold $M$ equipped with a general
connection $\nabla $ and with a metric tensor 
\slg%
. In Section 2 we introduced the notion of a parallelism structure on $%
U\subset M$, i.e., a pair $(U,\gamma )$ where $\gamma $ is a \emph{%
connection }extensor field on $U\subset M$ (representing there the
restriction in $U$ of some connection defined on $M$), and the theory of $a$%
-directional covariant derivatives of the representatives of multivector and
extensor fields on $U$, analyzing the main properties satisfied by these
objects. We also give in Section 2 a novel and intrinsic presentation (i.e.,
one that does not depend on a chosen orthonormal moving frame) of Cartan
theory of the torsion and curvature fields and of Cartan's structure
equations. We next introduce a metric metric structure $(U,g)$ in $U\subset
M $ and corresponding Christoffel operators\footnote{%
These objects generalize the Christofell symbols of the standard formalism
that are defined for coordinate basis vectors.} and the Levi-Civita
connection field. In Section 3 we studied in details the metric
compatibility of covariant derivatives and proved several important formulas
that are useful for applications of the formalism. In Section 4 several
operator derivatives in geometric structures $(U,\gamma ,g)$, like ordinary
and covariant Hodge coderivatives, are presented. The noticeable concept of
gauge (deformed) derivative is introduced and a crucial result (Theorem 2)
involving pairs of distinct metric compatible derivatives is proved . In
Section 4 we proved several useful duality identities. The relation between
all the concepts introduced has been carefully investigated.

\section{Appendix}

Let $U$ be an open subset of $U_{o},$ and let $(U,\phi)$ and $%
(U,\phi^{\prime })$ be two local coordinate systems on $U$ compatibles with $%
(U_{o},\phi _{o}).$ As we know there must be two pairs of reciprocal frame
fields on $U.$ The covariant and contravariant frame fields $%
\{b_{\alpha}\cdot\partial x_{o}\}$ and $\{\partial_{o}x^{\alpha}\}$
associated to $(U,\phi),$ and those ones $\{b_{\mu}\cdot\partial^{%
\prime}x_{o}\}$ and $\{\partial_{o}x^{\mu^{\prime}}\}$ associated to $%
(U,\phi^{\prime}).$

\textbf{A.1} The $n^{3}$ smooth scalar fields on $U,$ namely $\Gamma
_{\alpha\beta}^{\gamma},$ defined by 
\begin{equation}
\Gamma_{\alpha\beta}^{\gamma}=\Gamma^{+}(b_{\alpha}\cdot\partial
x_{o},b_{\beta}\cdot\partial x_{o})\cdot\partial_{o}x^{\gamma}  \tag{A1}
\end{equation}
correspond to the classically so-called \emph{coefficients of connection,}
of course, associated to $(U,\phi)$. The coefficients of connection
associated to $(U,\phi^{\prime})$ are given by 
\begin{equation}
\Gamma_{\mu^{\prime}\nu^{\prime}}^{\lambda^{\prime}}=\Gamma^{+}(b_{\mu}\cdot%
\partial^{\prime}x_{o},b_{\nu}\cdot\partial^{\prime}x_{o})\cdot
\partial_{o}x^{\lambda^{\prime}}.  \tag{A2}
\end{equation}

We check next what is the \emph{law of transformation} between them. We will
employ the simplified notations: $b_{\mu}\cdot\partial^{\prime}x_{o}=\dfrac{%
\partial x_{o}^{\sigma}}{\partial x^{\mu^{\prime}}}b_{\sigma}$ and $%
\partial_{o}x^{\alpha}=b^{\tau}\dfrac{\partial x^{\alpha}}{\partial
x_{o}^{\tau}},$ etc. Then, by recalling the expansion formulas for smooth
vector fields, $v=(v\cdot\partial_{o}x^{\alpha})b_{\alpha}\cdot\partial
x_{o} $ and $v=(v\cdot b_{\gamma}\cdot\partial x_{o})\partial_{o}x^{\gamma}$%
, using Eqs.(\ref{CO.2a}), (\ref{CO.2b}), (\ref{CO.2c}) and (\ref{CO.2d}),
and recalling the remarkable identity $(b_{\alpha}\cdot\partial
x_{o})\cdot\partial_{o}X=b_{\alpha}\cdot\partial X$, where $X\in\mathcal{M}%
(U)$ (see \cite{1}), we can write 
\begin{align*}
\Gamma_{\mu^{\prime}\nu^{\prime}}^{\lambda^{\prime}} & =\Gamma^{+}(\frac{%
\partial x^{\alpha}}{\partial x^{\mu^{\prime}}}b_{\alpha}\cdot\partial x_{o},%
\frac{\partial x^{\beta}}{\partial x^{\nu^{\prime}}}b_{\beta}\cdot\partial
x_{o})\cdot\frac{\partial x^{\lambda^{\prime}}}{\partial x^{\gamma}}%
\partial_{o}x^{\gamma} \\
& =\frac{\partial x^{\alpha}}{\partial x^{\mu^{\prime}}}\frac{\partial
x^{\beta}}{\partial x^{\nu^{\prime}}}\frac{\partial x^{\lambda^{\prime}}}{%
\partial x^{\gamma}}\Gamma_{\alpha\beta}^{\gamma} \\
& +\frac{\partial x^{\alpha}}{\partial x^{\mu^{\prime}}}(b_{\alpha}\cdot%
\partial x_{o})\cdot\partial_{o}(\frac{\partial x^{\beta}}{\partial
x^{\nu^{\prime}}})\frac{\partial x^{\lambda^{\prime}}}{\partial x^{\gamma}}%
(b_{\beta}\cdot\partial x_{o}\cdot\partial_{o}x^{\gamma}) \\
& =\frac{\partial x^{\alpha}}{\partial x^{\mu^{\prime}}}\frac{\partial
x^{\beta}}{\partial x^{\nu^{\prime}}}\frac{\partial x^{\lambda^{\prime}}}{%
\partial x^{\gamma}}\Gamma_{\alpha\beta}^{\gamma}+\frac{\partial x^{\alpha}}{%
\partial x^{\mu^{\prime}}}\frac{\partial}{\partial x^{\alpha}}(\frac{%
\partial x^{\beta}}{\partial x^{\nu^{\prime}}})\frac{\partial
x^{\lambda^{\prime}}}{\partial x^{\gamma}}\frac{\partial x^{\gamma}}{%
\partial x^{\beta}},
\end{align*}
i.e., 
\begin{equation}
\Gamma_{\mu^{\prime}\nu^{\prime}}^{\lambda^{\prime}}=\frac{\partial
x^{\alpha }}{\partial x^{\mu^{\prime}}}\frac{\partial x^{\beta}}{\partial
x^{\nu ^{\prime}}}\frac{\partial x^{\lambda^{\prime}}}{\partial x^{\gamma}}%
\Gamma_{\alpha\beta}^{\gamma}+\frac{\partial^{2}x^{\beta}}{\partial
x^{\mu^{\prime}}\partial x^{\nu^{\prime}}}\frac{\partial x^{\lambda^{\prime}}%
}{\partial x^{\beta}}.  \tag{A3}
\end{equation}
It is just the well-known law of transformation for the \emph{classical
coefficients of connection} associated to each of the \emph{coordinate
systems} \texttt{\ }$\{x^{\mu}\}$ and \texttt{\ }$\{x^{\mu^{\prime}}\}.%
\vspace{0.1in}$

\textbf{A.2} Given a smooth vector field on $U,$ say $v,$ the covariant and
contravariant components of $v$ with respect to $(U,\phi)$ and $(U,\phi
^{\prime})$ are respectively given by 
\begin{align}
v_{\alpha} & =v\cdot(b_{\alpha}\cdot\partial x_{o}),  \tag{A4} \\
v^{\alpha} & =v\cdot\partial_{o}x^{\alpha},  \tag{A5} \\
v_{\alpha}^{\prime} & =v\cdot(b_{\alpha}\cdot\partial^{\prime}x_{o}), 
\tag{A6} \\
v^{\alpha^{\prime}} & =v\cdot\partial_{o}x^{\alpha^{\prime}}.  \tag{A7}
\end{align}

We check the relationship between the covariant components of $v$ with
respect to each of the coordinate systems \texttt{\ }$\{x^{\mu}\}$ and 
\texttt{\ }$\{x^{\mu^{\prime}}\}.$ By recalling the expansion formula $%
w=(w\cdot \partial_{o}x^{\beta})b_{\beta}\cdot\partial x_{o},$ we have 
\begin{align}
v_{\alpha}^{\prime} &
=v\cdot(b_{\alpha}\cdot\partial^{\prime}x_{o}\cdot\partial_{o}x^{\beta})b_{%
\beta}\cdot\partial x_{o},  \notag \\
v_{\alpha}^{\prime} & =\frac{\partial x^{\beta}}{\partial x^{\alpha^{\prime
}}}v_{\beta}.  \tag{A8}
\end{align}
It is the expected law of transformation for the covariant components of $v$
associated to each of \texttt{\ }$\{x^{\mu}\}$ and \texttt{\ }$%
\{x^{\mu^{\prime }}\}.$

Analogously, by using the expansion formula $w=(w\cdot
b_{\beta}\cdot\partial x_{o})\cdot\partial_{o}x^{\beta},$ we can get the
classical law of transformation for the contravariant components of $v,$
i.e., 
\begin{equation}
v^{\alpha^{\prime}}=\frac{\partial x^{\alpha^{\prime}}}{\partial x^{\beta}}%
v^{\beta}.  \tag{A9}
\end{equation}
\vspace{0.1in}

\textbf{A.3} We see next which is the meaning of the classical covariant
derivatives of the contravariant and covariant components of a smooth vector
field. By using the expansion formula $v=v^{\alpha}b_{\alpha}\cdot\partial
x_{o},$ and Eqs.(\ref{CDM.4b}) and (\ref{CDM.4c}), and recalling once again
the remarkable identity $(b_{\mu}\cdot\partial
x_{o})\cdot\partial_{o}X=b_{\mu}\cdot\partial X,$ where $X\in\mathcal{M}(U),$
we can write 
\begin{align*}
(\nabla_{b_{\mu}\cdot\partial x_{o}}^{+}v)\cdot\partial_{o}x^{\lambda} &
=(b_{\mu}\cdot\partial
x_{o}\cdot\partial_{o}v^{\alpha})(b_{\alpha}\cdot\partial
x_{o}\cdot\partial_{o}x^{\lambda}) \\
& +v^{\alpha}(\nabla_{b_{\mu}\cdot\partial x_{o}}^{+}b_{\alpha}\cdot\partial
x_{o})\cdot\partial_{o}x^{\lambda} \\
& =\frac{\partial v^{\alpha}}{\partial x^{\mu}}\delta_{\alpha}^{\lambda
}+v^{\alpha}\Gamma^{+}(b_{\mu}\cdot\partial x_{o},b_{\alpha}\cdot\partial
x_{o})\cdot\partial_{o}x^{\lambda} \\
& =\frac{\partial v^{\lambda}}{\partial x^{\mu}}+\Gamma_{\mu\alpha}^{\lambda
}v^{\alpha},
\end{align*}
i.e., 
\begin{equation}
(\nabla_{b_{\mu}\cdot\partial
x_{o}}^{+}v)\cdot\partial_{o}x^{\lambda}=\left. v^{\lambda}\right. _{;\text{ 
}\mu}.  \tag{A10}
\end{equation}
By using the expansion formula $v=v_{\alpha}\partial_{o}x^{\alpha},$ Eqs.(%
\ref{CDM.4b}) and (\ref{CDM.4c}), and Eq.(\ref{CO.3}), we get 
\begin{equation*}
(\nabla_{b_{\mu}\cdot\partial x_{o}}^{-}v)\cdot b_{\nu}\cdot\partial x_{o}=%
\frac{\partial v_{\nu}}{\partial x^{\mu}}-\Gamma_{\mu\nu}^{\alpha
}v_{\alpha},
\end{equation*}
i.e., 
\begin{equation}
(\nabla_{b_{\mu}\cdot\partial x_{o}}^{-}v)\cdot b_{\nu}\cdot\partial
x_{o}=\left. v_{\nu}\right. _{;\text{ }\mu}.  \tag{A11}
\end{equation}
\vspace{0.1in}

\textbf{A.4} Now, for instance, let us take a smooth $(1,1)$-extensor field
on $U,$ say $t.$ The covariant and contravariant components of $t\ $with
respect to $(U,\phi)$ and $(U,\phi^{\prime})$ are respectively defined to be 
\begin{align}
t_{\mu\nu} & =t(b_{\mu}\cdot\partial x_{o})\cdot b_{\nu}\cdot\partial x_{o},
\tag{A12} \\
t^{\mu\nu} & =t(\partial_{o}x^{\mu})\cdot\partial_{o}x^{\nu},  \tag{A13} \\
t_{\mu^{\prime}\nu^{\prime}} & =t(b_{\mu}\cdot\partial^{\prime}x_{o})\cdot
b_{\nu}\cdot\partial^{\prime}x_{o},  \tag{A14} \\
t^{\mu^{\prime}\nu^{\prime}} & =t(\partial_{o}x^{\mu^{\prime}})\cdot
\partial_{o}x^{\nu^{\prime}}  \tag{A15}
\end{align}
It is also possible to introduce two mixed components of $t$ with respect to 
$(U,\phi)$ and $(U,\phi^{\prime}).$ They are defined by 
\begin{align}
\left. t_{\mu}\right. ^{\nu} & =t(b_{\mu}\cdot\partial x_{o})\cdot
\partial_{o}x^{\nu},  \tag{A16} \\
\left. t^{\mu}\right. _{\nu} & =t(\partial_{o}x^{\mu})\cdot
b_{\nu}\cdot\partial x_{o},  \tag{A17} \\
\left. t_{\mu^{\prime}}\right. ^{\nu^{\prime}} & =t(b_{\mu}\cdot
\partial^{\prime}x_{o})\cdot\partial_{o}x^{\nu^{\prime}},  \tag{A18} \\
\left. t^{\mu^{\prime}}\right. _{\nu^{\prime}} &
=t(\partial_{o}x^{\mu^{\prime}})\cdot b_{\nu}\cdot\partial^{\prime}x_{o}. 
\tag{A19}
\end{align}

We will try to check the law of transformation for the covariant component
of $t.$ We can write 
\begin{align}
t_{\mu^{\prime}\nu^{\prime}} & =t(b_{\mu}\cdot\partial^{\prime}x_{o})\cdot
b_{\nu}\cdot\partial^{\prime}x_{o}  \notag \\
& =t(\frac{\partial x^{\alpha}}{\partial x^{\mu^{\prime}}}%
b_{\alpha}\cdot\partial x_{o})\cdot\frac{\partial x^{\beta}}{\partial
x^{\nu^{\prime}}}b_{\beta}\cdot\partial x_{o}  \notag \\
& =\frac{\partial x^{\alpha}}{\partial x^{\mu^{\prime}}}\frac{\partial
x^{\beta}}{\partial x^{\nu^{\prime}}}t(b_{\alpha}\cdot\partial x_{o})\cdot
b_{\beta}\cdot\partial x_{o},  \notag \\
t_{\mu^{\prime}\nu^{\prime}} & =\frac{\partial x^{\alpha}}{\partial
x^{\mu^{\prime}}}\frac{\partial x^{\beta}}{\partial x^{\nu^{\prime}}}%
t_{\alpha\beta}.  \tag{A20}
\end{align}
This is the classical law of transformation for the covariant components of
a smooth $2$-tensor field from \texttt{\ }$\{x^{\mu}\}$ to \texttt{\ }$%
\{x^{\mu^{\prime}}\}.$

By following similar steps we can get the laws of transformation for the
contravariant and mixed components of $t.$ We have 
\begin{align}
t^{\mu ^{\prime }\nu ^{\prime }}& =\frac{\partial x^{\mu ^{\prime }}}{%
\partial x^{\alpha }}\frac{\partial x^{\nu ^{\prime }}}{\partial x^{\beta }}%
t^{\alpha \beta },  \tag{A21} \\
\left. t_{\mu ^{\prime }}\right. ^{\nu ^{\prime }}& =\frac{\partial
x^{\alpha }}{\partial x^{\mu ^{\prime }}}\frac{\partial x^{\nu ^{\prime }}}{%
\partial x^{\beta }}\left. t_{\alpha }\right. ^{\beta },  \tag{A22} \\
\left. t^{\mu ^{\prime }}\right. _{\nu ^{\prime }}& =\frac{\partial x^{\mu
^{\prime }}}{\partial x^{\alpha }}\frac{\partial x^{\beta }}{\partial x^{\nu
^{\prime }}}\left. t^{\alpha }\right. _{\beta }.  \tag{A23}
\end{align}%
They perfectly agree with the classical laws of transformation for the
respective contravariant and mixed components of a smooth $2$-tensor field
from $\{x^{\mu }\}$ to $\{x^{\mu ^{\prime }}\}$.\vspace{0.1in}

\textbf{A.5} We now obtain the relation between the covariant derivatives of 
$t$ and the classical concepts of covariant derivatives of the covariant,
contravariant and mixed components of a smooth $2$-tensor field. For
instance, we have 
\begin{align*}
& (\nabla_{b_{\mu}\cdot\partial x_{o}}^{++}t)(b_{\alpha}\cdot\partial
x_{o})\cdot b_{\beta}\cdot\partial x_{o} \\
& =b_{\mu}\cdot\partial x_{o}\cdot\partial_{o}(t(b_{\alpha}\cdot\partial
x_{o})\cdot b_{\beta}\cdot\partial x_{o})-t(\nabla_{b_{\mu}\cdot\partial
x_{o}}^{+}b_{\alpha}\cdot\partial x_{o})\cdot b_{\beta}\cdot\partial x_{o} \\
& -t(b_{\alpha}\cdot\partial x_{o})\cdot\nabla_{b_{\mu}\cdot\partial
x_{o}}^{+}b_{\beta}\cdot\partial x_{o} \\
& =\frac{\partial t_{\alpha\beta}}{\partial x^{\mu}}-t(\Gamma^{+}(b_{\mu
}\cdot\partial x_{o},b_{\alpha}\cdot\partial
x_{o})\cdot\partial_{o}x^{\sigma }b_{\sigma}\cdot\partial x_{o})\cdot
b_{\beta}\cdot\partial x_{o} \\
& -t(b_{\alpha}\cdot\partial x_{o})\cdot\Gamma^{+}(b_{\mu}\cdot\partial
x_{o},b_{\beta}\cdot\partial
x_{o})\cdot\partial_{o}x^{\tau}b_{\tau}\cdot\partial x_{o} \\
& =\frac{\partial t_{\alpha\beta}}{\partial x^{\mu}}-\Gamma_{\mu\alpha
}^{\sigma}t_{\sigma\beta}-\Gamma_{\mu\beta}^{\tau}t_{\alpha\tau},
\end{align*}
i.e., 
\begin{equation}
(\nabla_{b_{\mu}\cdot\partial x_{o}}^{++}t)(b_{\alpha}\cdot\partial
x_{o})\cdot b_{\beta}\cdot\partial x_{o}=\left. t_{\alpha\beta}\right. _{;%
\text{ }\mu}.  \tag{A24}
\end{equation}
We can also write 
\begin{align*}
& (\nabla_{b_{\mu}\cdot\partial x_{o}}^{+-}t)(b_{\alpha}\cdot\partial
x_{o})\cdot\partial_{o}x^{\beta} \\
& =b_{\mu}\cdot\partial x_{o}\cdot\partial_{o}(t(b_{\alpha}\cdot\partial
x_{o})\cdot\partial_{o}x^{\beta})-t(\nabla_{b_{\mu}\cdot\partial
x_{o}}^{+}b_{\alpha}\cdot\partial x_{o})\cdot\partial_{o}x^{\beta} \\
& -t(b_{\alpha}\cdot\partial x_{o})\cdot\nabla_{b_{\mu}\cdot\partial
x_{o}}^{-}\partial_{o}x^{\beta} \\
& =\frac{\partial\left. t_{\alpha}\right. ^{\beta}}{\partial x^{\mu}}%
-t(\Gamma^{+}(b_{\mu}\cdot\partial x_{o},b_{\alpha}\cdot\partial
x_{o})\cdot\partial_{o}x^{\sigma}b_{\sigma}\cdot\partial x_{o})\cdot\partial
_{o}x^{\beta} \\
& -t(b_{\alpha}\cdot\partial x_{o})\cdot\Gamma^{-}(b_{\mu}\cdot\partial
x_{o},\partial_{o}x^{\beta})\cdot b_{\tau}\cdot\partial x_{o}\partial
_{o}x^{\tau}. \\
& =\frac{\partial\left. t_{\alpha}\right. ^{\beta}}{\partial x^{\mu}}%
-\Gamma_{\mu\alpha}^{\sigma}\left. t_{\sigma}\right. ^{\beta}+\Gamma
_{\mu\tau}^{\beta}\left. t_{\alpha}\right. ^{\tau},
\end{align*}
i.e., 
\begin{equation}
(\nabla_{b_{\mu}\cdot\partial x_{o}}^{+-}t)(b_{\alpha}\cdot\partial
x_{o})\cdot\partial_{o}x^{\beta}=\left. \left. t_{\alpha}\right. ^{\beta
}\right. _{;\text{ }\mu}.  \tag{A25}
\end{equation}
\medskip

\textbf{Acknowledgments: } V. V. Fern\'{a}ndez and A. M. Moya are very
grateful to Mrs. Rosa I. Fern\'{a}ndez who gave to them material and
spiritual support at the starting time of their research work. This paper
could not have been written without her inestimable help. Authors are also
grateful to Drs. E. Notte-Cuello and E. Capelas de Oliveira for useful
discussions.


\begin{thebibliography}{9}
\bibitem{1} {\footnotesize Moya, A. M., Fern\'{a}ndez, V. V., and Rodrigues,
W. A., Jr. Multivector and Extensor Fields on Smooth Manifolds, \textit{Int.
J. Geom. Meth. Mod. Phys. \ }\textbf{4 }(6) (2007).}

\bibitem{2} {\footnotesize Fern\'{a}ndez, V. V., Moya, A. M., and Rodrigues,
W. A., Jr., \ Geometric Algebras and Extensors, \textit{Int. J. Geom. Meth.
Mod. Phys}. \ \textbf{4 }(6}) {\footnotesize (2007).}

\bibitem{3} {\footnotesize Fern\'{a}ndez, V. V., Moya, A. M., and Rodrigues,
W. A., Jr., Lagrangian Formalism for Multiform Fields on Minkowski
Spacetime, \textit{Int. J. Theor. Phys}. \textbf{40}, 299-313 (2001).}

\bibitem{rodoliv2006} {\footnotesize Rodrigues, W. A. Jr., and Oliveira, E.
Capelas, \textit{The Many Faces of Maxwell, Dirac and Einstein Equations},
to appear in Lecture Notes in Physics \textbf{722}, Springer, New York, 2007.%
}
\end{thebibliography}
\end{document}